\documentclass[12pt,letterpaper]{article}
\usepackage{amsfonts}
\usepackage{amssymb}
\usepackage{latexsym}  
\usepackage{pifont}  
\usepackage{amsmath}   

\DeclareMathAlphabet{\mathbbold}{U}{bbold}{m}{n}

\let\SavedRightarrow=\Rightarrow
\usepackage{marvosym}
\let\Rightarrow=\SavedRightarrow

\input cyracc.def
\newfam\cyrfam

\font\twelvecyr=wncyr10 scaled 1200 
\def\cyr{\fam\cyrfam\twelvecyr\cyracc}

\newcommand\Sh{\mbox{\cyr Sh}}  

\renewcommand\AA{{\mathcal A}}
\newcommand\BB{{\mathcal B}}
\newcommand\DD{{\mathcal D}}
\newcommand\UU{{\mathcal U}}
\newcommand\PP{{\mathcal P}} 

\newcommand\RRR{{\mathbb R}}
\newcommand\TTT{{\mathbb T}}
\newcommand\NNN{{\mathbb N}}
\newcommand\QQQ{{\mathbb Q}}
\newcommand\CCC{{\mathbb C}}
\newcommand\PPP{{\mathbb P}}

\newcommand\cccc{{\mathfrak c}}
\newcommand\KKKK{{\mathfrak K}}
\newcommand\RRRR{{\mathfrak R}}
\newcommand\MMMM{{\mathfrak M}}
\newcommand\MMMMS{{\mathfrak M}_{\mathfrak S}}

\newcommand\yyy{\mathbbold{y}}
\newcommand\sss{\mathbbold{s}}

\newcommand\cchi{{\raise 2 pt \hbox{$\chi$}}}

\newcommand\res{\mathord {\upharpoonright}}  

\newcommand\sdup{\ding{"4F}}

\newcommand\mop{\mathord{+}}  
\newcommand\mom{\mathord{-}} 

\newcommand\onto{\twoheadrightarrow}  

\newcommand\ran{\mathrm{ran}}  
\newcommand\bx{\mathrm{box}}  
\newcommand\lin{\mathrm{line}}  
\newcommand\dom{\mathrm{dom}}  
\newcommand\diam{\mathrm{diam}}   
\newcommand\supt{\mathrm{supt}}   
\newcommand\cl{\mathrm{cl}}   
\newcommand\lh{\mathrm{lh}}   
\newcommand\corn{\mathrm{corn}}  

\newcommand\ii{\mbox{\sc{i}}}
\newcommand\rr{\mbox{\sc{r}}}

\newcommand\lea{\sqsubseteq}  
\newcommand\leac{\, \mbox{\rlap{$\sqsubseteq$}  %
   \raise 2.3pt \hbox{\scriptsize $\mathrm{c}$}}\; }

\newcommand\cat{^{\mathord{\frown}}}  

\newcommand\iv{^{-1}} 

\def\eop{{\Large \Coffeecup}}  

\newenvironment{itemizz}{\begin{itemize}\setlength{\itemsep}{-1mm}} %
{\end{itemize}}                              

\newenvironment{itemizn}[1] 
{\begin{itemize} \setlength{\itemsep}{-1mm} %
\renewcommand\labelitemi{\ding{#1}}} %
{\end{itemize}}

\newtheorem{theorem}{Theorem}[section]
\newtheorem{definition}[theorem]{Definition}
\newtheorem{lemma}[theorem]{Lemma}

\newtheorem{corollary}[theorem]{Corollary}
\newtheorem{proposition}[theorem]{Proposition}
\newtheorem{question}[theorem]{Question}
\newtheorem{example}[theorem]{Example}

\newenvironment{proof}{{\bf Proof.}}{\eop\medskip}
\newenvironment{proofof}[1]{\medskip \textbf{Proof of #1.}}{\eop\medskip}

\begin{document}

\title{Ordered Spaces, Metric Preimages, and Function Algebras
}

\author{
Kenneth Kunen\footnote{University of Wisconsin,  Madison, WI  53706, U.S.A.,
\ \ kunen@math.wisc.edu}
\thanks{Author partially supported by NSF Grant DMS-0456653.}
}

\maketitle

\begin{abstract}
We consider the
\textit{Complex Stone-Weierstrass Property} (CSWP), which is
the complex version of the Stone-Weier\-strass Theorem.
If $X$ is a compact subspace of a product of three linearly ordered
spaces, then
$X$ has the CSWP if and only if
$X$ has no subspace homeomorphic to the Cantor set.
In addition, every finite power of the double arrow space
has the CSWP.  These results are proved using some results about those
compact Hausdorff spaces which have scattered-to-one
maps onto compact metric spaces.
\end{abstract}

\section{Introduction} 
\label{sec-intro}
All topologies discussed in this paper are assumed to be Hausdorff.
As usual, a subset of a space is \emph{perfect} iff it is closed and
non-empty and has no isolated points, so $X$ is \emph{scattered} iff
$X$ has no perfect subsets.

The usual version of the
\emph{Stone-Weier\-strass Theorem} involves
subalgebras of $C(X,\RRR)$, and is true for all compact $X$.
If one replaces the real numbers $\RRR$ by the complex numbers $\CCC$,
the ``theorem'' is true for some $X$ and false for others,
so it becomes a \emph{property} of $X$:

\begin{definition}
\label{def-cswp}
If $X$ is compact, then
$C(X) = C(X,\CCC)$ is the algebra of continuous
complex-valued functions on $X$,
with the usual supremum norm.
$\AA \lea C(X)$ means that $\AA$ is a subalgebra of $C(X)$
which separates points and contains the constant functions.
$\AA \leac C(X)$ means that $\AA \lea C(X)$ and
$\AA$ is closed in $C(X)$.
$X$ has the \emph{Complex Stone-Weierstrass Property (CSWP)}
iff every $\AA \lea C(X)$ is dense in $C(X)$;
equivalently,
iff every $\AA \leac C(X)$ equals $C(X)$.
\end{definition}

The CSWP is easily seen to be true for finite spaces.
The complex analysis developed in the 1800s shows that the CSWP is false
for many compact subspaces of the plane;
for example, it is false for the unit circle $\TTT$;
the classic counter-example being the algebra
of complex polynomials $\PP \lea C(\TTT)$.
These remarks are subsumed by results of W. Rudin \cite{RUD1,RUD2}
from the 1950s:

\begin{theorem}
\label{thm-rudin}
Let $X$ be any compact space.
\begin{itemizz}
\item[1.] If $X$ contains a copy of the
Cantor set, then $X$ fails the CSWP.
\item[2.] If $X$ is scattered, then $X$
satisfies the CSWP.
\end{itemizz}
\end{theorem}

If a compact space is metrizable (equivalently, second countable),
then it contains a Cantor subset iff it is not scattered,
so as Rudin pointed out:

\begin{corollary}
\label{cor-metric}
If $X$ is compact metric, then $X$ satisfies the CSWP
iff $X$ does not contain a copy of the Cantor set.
\end{corollary}

One might conjecture that this corollary holds for all
compact $X$, but that was refuted in 1960 by
Hoffman and Singer \cite{HS} (see also \cite{GA,HO});
their results imply that any compactum containing $\beta\NNN$
fails the CSWP.

However, the corollary does hold for some more ``reasonable''
classes of spaces.   Kunen \cite{KU} showed in 2004:

\begin{theorem}
\label{thm-lots}
If $X$ is a compact LOTS, then $X$ satisfies the CSWP
iff $X$ does not contain a copy of the Cantor set.
\end{theorem}

As usual, a LOTS is a linearly ordered topological space.
Of course, the $\rightarrow$
of this result is clear from Theorem \ref{thm-rudin};
only the $\leftarrow$ was new.
This theorem shows that there are some non-scattered 
spaces with the CSWP, such as 
the double arrow space of Alexandroff and Urysohn
(see Definition \ref{def-das}, or \cite{AU}, p.~76).

One can now ask whether there are further classes of ``\emph{reasonable}''
spaces for which results such as Corollary \ref{cor-metric}
and Theorem \ref{thm-lots} hold.  We do not know the best possible
result along this line, but we shall prove in Section \ref{sec-cswp}:

\begin{theorem}
\label{thm-three-lots}
If $X$ is compact and $X \subseteq L_0 \times L_1 \times L_2$,
where $L_0,L_1,L_2$ are LOTSes, then $X$ has the CSWP
iff $X$ does not contain a copy of the Cantor set.
\end{theorem}

Here, we may assume that $L_0,L_1,L_2$ are compact
(otherwise, replace them by the projections of $X$).
It is unknown whether the product of two spaces with the CSWP
must also have the CSWP.  Even if this turns out to be true,
Theorem \ref{thm-three-lots} is not immediate from Theorem \ref{thm-lots},
since $X$ is an arbitrary compact subset of the product,
and $L_0,L_1,L_2$ may fail the CSWP (i.e., have Cantor subsets).

By a slightly different argument, we shall show in Section \ref{sec-powers}:

\begin{theorem}
\label{thm-das-power}
If $L$ is the double arrow space, then $L^n$ has the CSWP
for every finite $n$.
\end{theorem}

Theorems \ref{thm-three-lots}  and \ref{thm-das-power} are proved
using some results from Section \ref{sec-tight}
about spaces which have
scattered-to-one maps onto metric spaces.
In Theorem \ref{thm-das-power}, there is a natural
$f : L^n \onto [0,1]^n$ for which the inverse of each point
is scattered (and of size $2^n$).
In Theorem \ref{thm-three-lots}, the $L_j$ need not have any
scattered-to-one maps onto metric spaces, but
a standard argument using measures reduces the proof of
Theorem \ref{thm-three-lots} to the case where the $L_j$ are separable
(see Section \ref{sec-cswp-reduce}), in which case
$X$ must have an eight-to-one map onto a compact
metric space.

If $L_0,L_1,L_2$ are separable in Theorem \ref{thm-three-lots},
then $X$ must also be first countable, and hence ``small'' in the cardinal
functions sense (see Juh\'asz \cite{JU}).
However, we do not believe that there is
a notion of ``reasonable'' involving only cardinal functions.
In \cite{HK2} it is shown that in some models of set theory,
there is a  compact $X$ which does not contain Cantor
subsets and which fails the CSWP, such that $X$ is both 
hereditarily separable and hereditarily Lindel\"of (and hence
also first countable).  In these models, $2^{\aleph_0} = \aleph_1$
and the standard cardinal functions of our $X$
(all either $\aleph_0$ or $\aleph_1$) are the least possible
among non-metric compacta.

Section \ref{sec-LOTS} reviews some elementary fact about LOTSes.
Section \ref{sec-remov} discusses the notion of a
\emph{removable space} defined in \cite{HK1}; this is a strengthening
of the CSWP used in Section \ref{sec-powers}.

\begin{definition} 
\label{def-local}
Let $\KKKK$ be a class of compact spaces.
$\KKKK$ is \emph{closed-hereditary} iff every closed
subspace of a space in $\KKKK$ is also in $\KKKK$.
$\KKKK$ is \emph{local} iff $\KKKK$ is closed-hereditary \emph{and}
for every compact $X$:  if $X$ is covered by open sets whose
closures lie in $\KKKK$, then $X \in \KKKK$.
\end{definition} 

Classes of compacta which restrict cardinal functions
(first countable, second countable, countable tightness, etc.)
are clearly local, whereas the class of compacta which
are homeomorphic to a LOTS is closed-hereditary, but not local.

It is easily seen that the CSWP is closed-hereditary;
this is Lemma 1.3 of \cite{KU}, but
the proof is implicit in Rudin \cite{RUD1}.
Thus, to prove part (1) of Theorem \ref{thm-rudin} in \cite{RUD1},
it was sufficient to show that the Cantor set itself fails the CSWP.

The removable spaces form a local class 
(see Section \ref{sec-remov}).
It is unknown whether the CSWP is a local property.
A proof that it is local cannot be completely trivial.
For example, locality would imply that the failure of the CSWP for $\TTT$
yields the failure of the CSWP for an arc $A \subseteq \TTT$.
Now, $A$ does in fact fail the CSWP, since it contains a Cantor set,
but we do not know how to construct a counter-example on $A$ directly
from the polynomial algebra $\PP \lea C(\TTT)$; note that the restriction
$\PP \res A \lea C(A)$ is dense in $C(A)$ by Mergelyan's Theorem.

\section{Ordered Spaces}
\label{sec-LOTS}
We begin by defining the double arrow space and some variants thereof:

\begin{definition}
\label{def-das}
$I = [0,1]$.  If $\Lambda : I \to \omega$, then
$I_\Lambda = \bigcup_{x \in I} \{x\} \times \{0,1,\ldots, \Lambda(x)\}$,
which is given the lexicographic order and the usual order topology.
If $S \subseteq (0,1)$, then $I_S  = I_{\cchi_S}$, where
$\cchi_S$ is the characteristic function;
then for $x \in S$, let $x^- = (x,0)$ and $x^+ = (x,1)$; 
while  if $x \notin S$, let $x^- = x^+ = (x,0)$.
The \emph{double arrow space} is $I_{(0,1)}$.
For any $\Lambda$, the map $(x, \ell) \mapsto x$ 
is the \emph{standard map} from $I_\Lambda$ onto $I$.
\end{definition}

So, we form $I_\Lambda$ by splitting each $x \in S$ into
$\Lambda(x) + 1$ neighboring points.  For $I_S$,
we split each $x \in S$ into two
neighboring points, $x^-,x^+$, and we don't split the points
in $I \backslash S$;
it is convenient to have $x^\pm$ defined for \emph{all} $x \in I$,
so, for example, we can say
that for all $a < b$ in $I$, $(a^+, b^-)$ is an open interval in $I_S$.
$I_S$ has no isolated points because
$0,1\notin S$.  The double arrow space is obtained by
splitting all points other than $0,1$.
$I_\emptyset \cong I$, and $I_{\QQQ \cap (0,1)}$ is homeomorphic
to the Cantor set.

\begin{lemma}
For each $S \subseteq (0,1)$, $I_S$ is a compact separable LOTS
with no isolated points.  $I_S$ is second countable iff $S$ is countable.
Every $I_\Lambda$ is a compact first countable LOTS.
\end{lemma}

$I_\Lambda$ will not be separable unless $\{x : \Lambda(x) > 1\}$
is countable.
The study of compact separable LOTSes can
be reduced to spaces of the form $I_S$.  First note,
by Lutzer and Bennett \cite{LB}:

\begin{lemma}
\label{lemma-hs}
If $X$ is a separable LOTS, then $X$ is hereditarily separable
and hereditarily Lindel\"of.
\end{lemma}

Also, it is easy to check:

\begin{lemma}
\label{lemma-rel}
If $X$ is a LOTS and $H$ is a compact
subset of $X$, then the relative topology and the order topology
agree on $H$.
\end{lemma}

Relating this to our $I_S$:

\begin{lemma}
\label{lemma-LOTS-standard}
Let $X$ be a compact separable LOTS.  Then
\begin{itemize}
\item[1.] If $X$ is perfect, then $X$ is
homeomorphic to $I_S$ for some $S \subseteq (0,1)$.
\item[2.] If $X$ is not second countable,
then $X$ has a closed subspace which is 
homeomorphic to $I_S$ for some uncountable $S \subseteq (0,1)$. 
\item[3.] $X$ is homeomorphic to
a subset of $I_S$ for some $S \subseteq (0,1)$.
\end{itemize}
\end{lemma}
\begin{proof}
For (1):  Let $E \subseteq X$ be countable and dense in $X$ and contain the
first and last elements of $X$.  Let $B$ be the set of all $b \in E$
such that for some $a \in E$: $a < b$ and $(a,b) = \emptyset$.
Let $D = E \backslash B$.  Since $X$ has no isolated points,
$D$ is also dense in $X$ and contains the
first and last elements of $X$, and is also densely ordered.
Let $f$ be an order isomorphism from $D$ onto $\QQQ \cap [0,1]$.
Then $f$ extends in a natural way to a continuous $F: X \onto [0,1]$,
and $1 \le |F\iv\{r\}| \le 2$ for each $r \in [0,1]$.
Let $S = \{r : |F\iv\{r\}| = 2\}$.

For (2): Since $X$ is hereditarily Lindel\"of, the
Cantor-Bendixson sequence of $X$ has countable length
and removes countably many points.  Thus, $X$ is not scattered,
and, letting $H$ be the perfect kernel of $X$, $X \backslash H$ is countable.
Then $H$ is separable and not second countable,
so $H \cong I_S$ for some uncountable $S$.

For (3): Apply (1) to the space obtained from $X$ by replacing each
isolated point by a copy of the double arrow space.
\end{proof}

Note that $(I_S)^2$ is separable, but it is not hereditarily separable
when $S$ is uncountable; in fact, more general $I_\Lambda$ occur
naturally in such products.
Fixing an uncountable $S \subseteq (0,1)$, let
$L_n = I_{\Lambda_n}$, where
$\Lambda_n(x) = n$ for $n \in S$, and 
$\Lambda_n(x) = 0$ for $n \notin S$.
Then $L_n$ is not separable whenever $n \ge 2$,
and the diagonal of $(I_S)^k$ is homeomorphic to $L_{2^k -1}$.

\section{Tight Maps and Dissipated Spaces}
\label{sec-tight}
We recall some definitions and results from \cite{KU2}.
As usual,
$f : X \to Y$ means that $f$ is a \emph{continuous} map from $X$ to $Y$,
and $f : X \onto Y$ means that $f$ is a continuous map from $X$
\emph{onto} $Y$.

\begin{definition}
\label{def-tight}
Assume that $X,Y$ are compact and $f: X \to Y$.
\begin{itemizn}{43}
\item A \emph{loose family} for $f$ is a disjoint family $\PP$ of closed
subsets of $X$ such that for some non-scattered $Q \subseteq Y$,
$Q = f(P)$ for all $P \in \PP$.
\item $f$ is $\kappa$--\emph{tight} iff there are no loose families 
for $f$ of size $\kappa$.
\item $f$ is \emph{tight} iff $f$ is $2$--tight.
\end{itemizn}
\end{definition}

This notion gets weaker as $\kappa$ gets bigger.
$f$ is $1$--tight iff $f(X)$ is scattered,
so that ``$2$--tight'' is the first non-trivial case.
$f$ is trivially $|X|^+$--tight. 
The usual projection from $[0,1]^2$ onto $[0,1]$ is not
$2^{\aleph_0}$--tight.

Some easy equivalents to ``$\kappa$--tight'' are described
in Lemma 2.2 of \cite{KU2}:

\begin{lemma}
\label{lemma-tight-equiv}
Assume that $X,Y$ are compact and $f : X \to Y$.  Then
$(1) \leftrightarrow (2)$.
If $\kappa$ is finite
and $Y$ is metric, then all four of the following are equivalent:
\begin{itemizz}
\item[1.] There is a loose family of size $\kappa$. 
\item[2.]
There is a disjoint family $\PP$ of perfect subsets of $X$ with $|\PP| = \kappa$
and a perfect $Q \subseteq Y$ such that $Q = f(P)$ for all $P \in \PP$.
\item[3.] For some metric $M$ and $\varphi\in C(X,M)$,
$\{y \in Y : |\varphi(f\iv\{y\})| \ge \kappa\}$ is uncountable.
\item[4.] Statement $(3)$, with $M = [0,1]$.
\end{itemizz}
\end{lemma}

If $X,Y$ are both compact metric, then $f: X \to Y$ is $\kappa$--tight iff
$\{y \in Y : |f\iv\{y\}| \ge \kappa\}$ is countable
(see Theorem 2.7 of \cite{KU2}).
Of course, the $\leftarrow$ direction is trivial.
The $\rightarrow$ direction for non-metric $X$ and $\kappa = 2$
is refuted by the standard map from the double arrow space onto $[0,1]$,
which is tight by Lemma 2.3 of \cite{KU2}:

\begin{lemma}
\label{lemma-tight-LOTS}
If $X,Y$ are compact LOTSes and $f:X \to Y$ is order-preserving
$(x_1 < x_2 \to f(x_1) \le f(x_2))$,
then $f$ is tight.
\end{lemma}

One can estimate the tightness of product maps using
Lemma 2.14 of \cite{KU2}:

\begin{lemma}
\label{lemma-tight-prod}
Assume that for $i = 0,1$:
$X_i, Y_i$ are compact,
$f_i:X_i \to Y_i$ is $(m_i + 1)$--tight,
$m_i \le n_i < \omega$, and
$|f_i\iv\{y\}| \le n_i$ for \emph{all} $y \in Y_i$.
Then $f_0 \times f_1 : X_0 \times X_1 \to Y_0 \times Y_1$ 
is $(\max(m_0 n_1 , m_1 n_0) + 1)$-- tight.
\end{lemma}

The notion of a \emph{dissipated} compactum
(Definition \ref{def-dis} below) 
involves tight maps onto metric compacta, ordered by
\emph{fineness}, so we define:

\begin{definition}
\label{def-fine}
Assume that $X,Y,Z$ are compact, $f : X \to Y$, and
$g : X \to Z$.  Then $f \le g$, or
$f$ is \emph{finer than} $g$, iff there is a
$\Gamma \in C(f(X),g(X))$ such that $g = \Gamma \circ f$.
\end{definition}

\begin{lemma}
\label{lemma-fine-equiv}
Assume that $X,Y,Z$ are compact, $f : X \to Y$, and
$g : X \to Z$.  Then $f \le g$ iff
$\; \forall x_1, x_2 \in X\, [ f(x_1) = f(x_2) \to g(x_1) = g(x_2)]$.
\end{lemma}

\begin{definition}
Assume $X$ is compact.  Let $\MMMM(X)$, the \emph{metric projections}
of $X$, be the class of all maps $\pi$ such that $\pi : X \to Y$ 
for some compact metric $Y$.
Then $\pi \in \MMMMS(X) \subseteq \MMMM(X)$ iff in addition,
each $\pi\iv\{y\}$ is scattered.
\end{definition}

\begin{lemma}
\label{lemma-MS-down}
If $\pi,\sigma\in\MMMM(X)$ and $\pi \le \sigma \in \MMMMS(X)$,
then $\pi \in \MMMMS(X)$.
\end{lemma}

Observe that in the definition of $f \le g$, it is irrelevant
whether $f,g$ map $X$ \emph{onto} $Y,Z$.
Here, and in the definition of $\MMMM(X)$, we should really regard
$f$ in the set-theoretic sense as a set of ordered pairs,
not as a triple $(f,X,Y)$, so that
$f: X \to Y$ and $f: X \onto f(X)$ are exactly the same object.
One could also define $\MMMM(X)$ and $\MMMMS(X)$ as sets of closed
equivalence relations on $X$.

\begin{lemma}
\label{lemma-count}
$\MMMM(X)$ is countably directed.  That is,
if $\sigma_n \in \MMMM(X)$ for $n \in \omega$, then there is
a $\pi \in \MMMM(X)$ with $\pi \le \sigma_n$ for each $n$.
\end{lemma}

\begin{lemma}
\label{lemma-sing}
If $\sigma \in \MMMMS(X)$, then there is a $\pi \in \MMMMS(X)$
with $\pi \le \sigma$ and $\pi : X \onto Y$,
such that $\pi\iv\{b\}$ is a singleton for some $b\in Y$.
\end{lemma}
\begin{proof}
Say $\sigma : X \onto Z$.
Fix any $c \in Z$, and then fix $a \in \sigma\iv\{c\}$ such that
$a$ is isolated in $\sigma\iv\{c\}$.  Since $Z$ is metric,
$\{a\}$ is a $G_\delta$ in $X$, so fix any $f\in C(X, [0,1])$
with $\{a\} = f\iv\{1\}$.
Choose $\pi\in\MMMMS(X)$ with $\pi \le \sigma$ and $\pi\le f$.
\end{proof}

Only a scattered compactum $X$ has the property that \emph{all} maps
in $\MMMM(X)$ are tight: If $X$ is not scattered, then $X$ maps onto 
$[0,1]^2$; if we follow that map by the usual projection onto $[0,1]$,
we get a map from $X$ onto $[0,1]$ which is not
even $\cccc$--tight.
The \emph{dissipated} compacta have the property that
\emph{cofinally} many of these maps are tight:

\begin{definition}
\label{def-dis}
$X$ is \emph{$\kappa$--dissipated} iff $X$ is compact and
whenever $g \in \MMMM(X)$, 
there is a finer $\kappa$--tight $f \in \MMMM(X)$.
$X$ is \emph{dissipated} iff $X$ is $2$--dissipated.
\end{definition}

So, the $1$--dissipated compacta are the scattered compacta.
Metric compacta are dissipated because we can let $f$ be identity map.
By Lemma 3.12 of \cite{KU2}:

\begin{lemma}
\label{lemma-local}
For any $\kappa$, the class of $\kappa$--dissipated compacta
is a local class.
\end{lemma}

An easy example of a dissipated space is given by:

\begin{lemma}
\label{lemma-lots-dis}  
If $X$ is a compact LOTS, then $X$ is dissipated
\end{lemma}

The proof (see Lemma 3.4 of \cite{KU2}) shows that given 
$g \in \MMMM(X)$,
there is a finer  $f \in \MMMM(X)$
such that $f(X)$ is a compact metric LOTS and $f$ is order-preserving.

Note that just having \emph{one} tight map $g$ from $X$ onto some
metric compactum $Z$ is not sufficient to prove that 
$X$ is dissipated, since the tightness of $g$ says nothing at all about
the complexity of a particular $g\iv\{z\}$. 
However, if all $g\iv\{z\}$ are scattered, then just one tight $g$ is
enough by Lemma 3.5 of \cite{KU2}:

\begin{lemma}
\label{lemma-one-enough}
Assume that some $g \in \MMMMS(X)$ is $\kappa$--tight.
Then all $f \le g$ are also  $\kappa$--tight,
so that $X$ is $\kappa$--dissipated.
\end{lemma}

This suggests the following definition:

\begin{definition}
\label{def-super}
$\pi\in\MMMM(X)$ is $\kappa$--\emph{supertight}
iff $\pi$ is $\kappa$--tight and $\pi\in\MMMMS(X)$.
Then $X$ is $\kappa$--\emph{superdissipated} iff
some $\pi\in\MMMMS(X)$ is $\kappa$--supertight.
\end{definition}

Using Lemmas \ref{lemma-one-enough},  \ref{lemma-local},
and \ref{lemma-MS-down} above:

\begin{lemma}
\label{lemma-lim-finer}
If $\pi, \sigma \in\MMMMS(X)$, 
$\pi \le \sigma$, and $\sigma$ is $\kappa$--supertight,
then $\pi$ is $\kappa$--supertight.
\end{lemma}

\begin{lemma}
A compactum $X$ is $\kappa$--superdissipated iff
$X$ is $\kappa$--dissipated and $\MMMMS(X) \ne \emptyset$.
\end{lemma}

\begin{lemma}
\label{lemma-clh}
The class of $\kappa$--superdissipated compacta is a local class.
\end{lemma}

By Lemma \ref{lemma-tight-LOTS}:

\begin{lemma}
The standard map
$\sigma : I_\Lambda \onto I$ is $2$--supertight.
\end{lemma}

The situation for products is more complicated.
By Lemma \ref{lemma-tight-prod} and induction:

\begin{lemma}
\label{lemma-big-lim}
For any $n \ge 1$ and
$S_i \subseteq I $ \textup(for $i < n$\textup)\textup:
The standard map $\sigma : \prod_{i <n} I_{S_i} \onto I^n$ is
$(2^{n-1} + 1)$--supertight.
\end{lemma}

This result is best possible by Theorem 3.9 of \cite{KU2};
a product $\prod_{i< n}X_i$ is not $(2^{n-1} )$--dissipated
if each $X_i$ is a compact separable LOTS,
none of the $X_i$ is scattered, and
at most one of the $X_i$ is second countable.

\begin{definition}
\label{def-kernel}
The \emph{perfect kernel}, $\ker(X)$, is $\emptyset$ if $X$ is scattered,
and the largest perfect subset of $X$ otherwise.
\end{definition}

By Lemma \ref{lemma-tight-equiv}, the tightness of $\pi : X \to Y$
can be expressed using perfect subsets of $X$, so that

\begin{lemma}
\label{lemma-ker-tight}
$\pi : X \to Y$ is
$\kappa$--\,\textup(super\textup)tight iff
$\pi\res\ker(X)$  is
$\kappa$--\,\textup(super\textup)tight, and
the space $X$ is
$\kappa$--\,\textup(super\textup)dissipated iff
$\ker(X)$ is
$\kappa$--\,\textup(super\textup)dissipated.
\end{lemma}

\begin{lemma}
\label{lemma-get-one}
Assume that $\pi: X \to Y$ is $(n+2)$--supertight,
where $n \in \omega$, $X$ is compact and $Y$ is compact metric,
and $\{P_0, \ldots, P_n\}$ is a loose family for $\pi$
of size $n+1$, with each $\pi(P_j) = Q$.  Then 
each  $\pi\res P_j : P_j \to Q$ is $2$--supertight,
$\ker (\pi\iv(Q)) \subseteq \bigcup_j P_j$,
and $\pi\iv(Q)$ is $2$--superdissipated.
\end{lemma}
\begin{proof}
If tightness fails for $\pi\res P_j$, then we could find uncountable
closed $Q' \subseteq Q$ and disjoint closed $P_j^0, P_j^1 \subseteq P_j$
with $\pi(P_j^0) = \pi(P_j^1) = Q'$.  If $P'_k = P_k \cap \pi\iv(Q')$,
then the sets $P'_0, \ldots, P'_{j-1}, P_j^0, P_j^1, P'_{j+1}, \ldots, P'_n$
would be a loose family for $\pi$ of size $n+2$.
If $\ker (\pi\iv(Q)) \not\subseteq \bigcup_j P_j$,
we could find a perfect $R \subseteq \pi\iv(Q) \setminus \bigcup_j P_j$;
then $\pi(R)$ is non-scattered (since all $\pi\iv\{y\}$ are scattered),
and $R$ plus the $P_j$ would contradict the $(n+2)$--tightness of $\pi$.
Finally, by Lemma \ref{lemma-ker-tight}, it is sufficient
to prove that $\bigcup_j P_j$ is superdissipated, and this is done using
the map into $Y \times \{0,1,\ldots,n\}$ which sends
$x \in P_j$ to $(\pi(x),j)$.
\end{proof}

Finally, we mention two lemmas for the case that 
$X$ does not contain a Cantor subset.
$\pi: X \to Y$ is trivially $n$--supertight when
all $|\pi\iv\{y\}| < n$, but also

\begin{lemma}
\label{lemma-lim-prop}
Assume that $\pi: X \to Y$,
$X$ is compact, $Y$ is compact metric, each $|\pi\iv\{y\}| \le n$,
and $X$ has no Cantor subsets.  Then
$\pi$ is $n$--supertight.
\end{lemma}
\begin{proof}
If not, let $P_0, \ldots, P_{n-1} \subseteq X $ be a loose family, with each
$\pi(P_j) = Q$.  Then $Q$ has a Cantor subset, and
each $P_j$ is homeomorphic (via $\pi$) to $Q$.
\end{proof}

By the next lemma, the spaces $X$ we consider are always totally disconnected:

\begin{lemma}
\label{lemma-discon}
Assume that $\pi : X \to Y$, where $X$ is compact and $Y$ is metric.
Assume that each $\pi\iv\{y\}$ is totally disconnected and
$X$ does not contain a copy of the Cantor set.
Then $X$ is totally disconnected.
\end{lemma}
\begin{proof}
Assume that $X$ is not totally disconnected.
Fix a metric on $Y$ for which $\diam(Y) \le 1$.
Obtain $K_s$ for $s \in 2^{<\omega}$ to satisfy:
\begin{itemizz}
\item[1.] $K_s$ is an infinite closed connected subset of $X$.
\item[2.] $\diam(\pi(K_s)) \le 2^{- \lh(s)}$.
\item[3.] $K_{s\cat 0}, K_{s\cat 1} \subset K_s$
and $K_{s\cat 0} \cap K_{s\cat 1} = \emptyset$.
\end{itemizz}
Assuming that
this can be done, define $K_f = \bigcap_{n \in \omega} K_{f\res n}$
for $f \in 2^\omega$.  By (2), $|\pi(K_f)| = 1$; say $\pi(K_f) = \{y_f\}$.
But $K_f$ is connected and $\pi\iv\{y_f\}$ is totally disconnected,
so $|K_f| = 1$; say $K_f = \{x_f\}$.
Then $f \mapsto x_f$ is a homeomorphism from $2^\omega$ into $X$,
contradicting our assumptions about $X$.

To build the $K_s$:  For $K_{(\,)}$, just use the assumption that
$X$ is not totally disconnected.
Now, say we are given $K_s$.
Choose $x_0,x_1\in K_s$ with $x_0 \ne x_1$.  Then find disjoint
relatively open $U_0,U_1 \subseteq K_s$ with each $x_\ell \in U_\ell$
and $\diam(\pi(U_\ell)) \le 2^{- \lh(s) - 1}$.  Then find relatively open
$V_\ell  \subseteq K_s$
with $x_\ell \in V_\ell \subseteq  \overline{V_\ell}\subseteq U_\ell$.
Then, let $K_{s\cat \ell}$ be the connected component of
the point $x_\ell$ in the space $\overline{V_\ell}$,
and note that $K_{s\cat \ell}$ cannot be a singleton.
\end{proof}

\section{The CSWP: Two Reductions} 

\label{sec-cswp-reduce}
These reductions were described in \cite{KU}:
Using the standard theory of function algebras
(see \cite{GAM, GA}),
we can reduce the CSWP to the study of idempotents, and
we can reduce the study of the CSWP in LOTSes to the separable case.

If $f \in C(X)$, then $f$ is an \textit{idempotent} iff $f^2 = f$;
equivalently iff $f$ is the characteristic
function of some clopen set.
An idempotent is called \emph{nontrivial} iff it is not the
identically $0$ or the identically $1$ function.
As with other proofs of the CSWP \cite{HK1, KU},
we shall proceed by considering idempotents.
Following \cite{KU},

\begin{definition} 
The compact space $X$
has the \emph{NTIP} 
iff every $\AA \leac C(X)$ contains a non-trivial idempotent.
\end{definition}

So, the NTIP is trivially false of connected spaces.
If $X$ is not connected, then the CSWP implies the NTIP\@.
The following is Lemma 3.5 of \cite{KU};
it is also easy to prove from the Bishop
Antisymmetric Decomposition
(see \cite{BIS}, or Theorem 13.1 in Chapter II of \cite{GAM}).

\begin{lemma}
\label{lemma-restrict-to-perfect}
Assume that $X$ is compact
and every perfect subset of $X$ has the NTIP\@.  Then $X$ has the CSWP.
\end{lemma}

Among the totally disconnected spaces, the NTIP is strictly weaker
than the CSWP (see \cite{KU}).  However, the lemma implies the
following corollary, which is used to reduce proofs of the CSWP
to proofs of the weaker NTIP:

\begin{corollary}
\label{cor-idem}
If $\KKKK$ is a closed-hereditary class of compact spaces
and every perfect space in $\KKKK$ has the NTIP,
then every space in $\KKKK$ has the CSWP.
\end{corollary}

In particular, if $\KKKK$ is the class of compact scattered spaces,
then this corollary applies vacuously, so all spaces in $\KKKK$
have the CSWP.  If $\KKKK$ contains some non-scattered spaces,
then, as in \cite{KU, HK1}, we produce idempotents using:

\begin{lemma}
\label{lemma-get-clopen}
Suppose that $\AA \leac C(X)$ and
there is some  $h \in \AA$ such that either $\Re(h(X))$ or $\Im(h(X))$
is not connected.
Then $\AA$ contains a non-trivial idempotent.
\end{lemma}

This is easy to prove using Runge's Theorem;
see Lemma 2.5 of \cite{KU}, but the method was also used in
\cite{RUD1} and \cite{HS}.

It remains to describe how to obtain such an $h$.
If $X$ is scattered, then $\Re(h(X))$ is scattered also,
so any $h$ for which $\Re(h(X))$ is not a singleton will do;
this is essentially the argument of \cite{RUD2}.
In some other cases, we can obtain $h$ using a tight map of $X$ onto
a metric space; this is described in Section \ref{sec-cswp}.

We now turn to the second reduction.
As in \S5 of \cite{KU},

\begin{definition}
If $\mu$ is a regular complex Borel measure on the compact space $X$,
then $|\mu|$ denotes its total variation, and
$\supt(\mu) = \supt(|\mu|)$ denotes its \textup(\!closed\textup) support;
that is,
$\supt(\mu) = X \setminus \bigcup \{U\subseteq X :\; $U$ \mbox{ is open }
\ \ \& \ \ |\mu|(U) = 0\}$.
\end{definition} 

Considering measure orthogonal to $\AA$, we get:

\begin{lemma}
\label{lemma-supt-cswp}
Assume that $X$ is compact and that $\supt(\mu)$ has the CSWP
for all regular Borel measure $\mu$.
Then $X$ has the CSWP.
\end{lemma}

By  Corollary 5.4 of \cite{KU2}, every such $\supt(\mu)$ is separable
in the case that $X$ is $\aleph_0$--dissipated; for a LOTS $X$,
this was a much earlier folklore result. 

\begin{corollary}
\label{cor-CSWP-sep}
If $X$ fails the CSWP and is $\aleph_0$--dissipated, then 
some compact separable subspace of $X$ fails the CSWP.
\end{corollary}

\begin{question}
Is there a compact space $X$ which fails the CSWP such that
all compact separable subspaces of $X$ satisfy the CSWP?
\end{question}

This $X$ cannot be one of the three examples already
known to fail the CSWP --- namely, any space containing either the
Cantor set \cite{RUD1} or $\beta\NNN$ \cite{HS} or 
the examples of \cite{HK2,HK3} (obtained assuming $\diamondsuit$ or CH),
since all these spaces are separable.

Now, considering products of LOTSes:

\begin{lemma}
\label{lemma-prod-cswp}
Assume that $X$ is a compact subset of $\prod_{\alpha < \kappa} L_\alpha$,
where each $L_\alpha$ is a LOTS, and assume
that $X$ does not have the CSWP.  Then for some separable closed compact
$H_\alpha \subseteq L_\alpha$, the space
$X \cap \prod_{\alpha < \kappa} H_\alpha$ also fails the CSWP.
\end{lemma}
\begin{proof}
Let $\pi_\alpha : X \to L_\alpha$ be the usual coordinate projection.
We may assume that each $L_\alpha = \pi_\alpha(X)$, so that
$L_\alpha$ is compact.  Fix $\mu$ on $X$ such that
$\supt(\mu)$ fails the CSWP,  let
$\mu \pi_\alpha\iv$ be the induced measure on $L_\alpha$,
let $H_\alpha = \supt(\mu \pi_\alpha\iv)$,
which is separable,
and note that $\supt(\mu) \subseteq X \cap \prod_{\alpha < \kappa} H_\alpha$.
\end{proof}

\begin{lemma}
\label{lemma-one-S}
For any $\kappa \le \omega$:  Suppose that there is a compact 
$X \subseteq \prod_{\alpha < \kappa} L_\alpha$,
where each $L_\alpha$ is a LOTS, $X$ has no Cantor subset, and
$X$ does not have the CSWP.   Then there is such an $X$ which is
a subset of $(I_S)^\kappa$ for some $S \subseteq (0,1)$.
\end{lemma}
\begin{proof}
By Lemma \ref{lemma-prod-cswp}, we may assume that each $L_\alpha$
is separable and compact.  Now, let $L$ be the compact
separable LOTS obtained by
placing the $L_\alpha$ end-to-end, adding a point $\infty$ in the
case that $\kappa  = \omega$.  Then we may assume that 
$X \subseteq L^\kappa$.  Finally, replace $L$ by an 
$I_S$ using Lemma \ref{lemma-LOTS-standard}(3).
\end{proof}

\section{The CSWP and Tightness} 
\label{sec-cswp}
We show here how one can use the concepts from Section
\ref{sec-tight} to produce
idempotents, and thus to prove the CSWP.

\begin{definition}
\label{def-phat}
Assume that $\pi: X \to Y$, where $X,Y$ are compact.
Then, for $f \in C(X)$, define
$\widehat f = (\pi\times f)(X)$; that is, 
\[
\widehat f \ =\  \{(\pi(x), f(x)) : x \in X \} 
\  \subseteq\  Y \times \CCC \ \ .
\]
\end{definition}

\begin{lemma}
Each $\widehat f$ is compact.
\end{lemma}

We plan to apply the next definition and lemma to
sets of the form $\widehat f$:

\begin{definition}
Fix $E \subseteq Y \times \CCC$ and
$\Phi: \CCC^m \to \CCC$.  Then $E_y = \{z : (y,z) \in E\}$ and
\[
\Phi * E \ =\  \bigcup_{y \in Y} \Phi( (E_y)^m ) \ \subseteq\  \CCC \ \ .
\]
\end{definition}

\begin{lemma}
\label{lemma-open-cover}
Suppose that $F \subseteq Y \times \CCC$ is compact and
$\Phi: \CCC^m \to \CCC$ is continuous.
Let $\BB$ be an open base for $Y \times \CCC$
which is closed under finite unions.  Then $\Phi * F$ is compact and
$\Phi * F = \bigcap\{ \Phi * \overline U : U \in \BB \ \&\ F \subseteq U \}$.
\end{lemma}

\begin{lemma}
\label{lemma-scat-phi}
Assume that $\pi: X \to Y$ is $n$--supertight and $f \in C(X)$.
Fix a continuous $\Phi : \CCC^{n} \to \CCC$ such that
$\Phi(z_1, \ldots, z_n) = 0$ unless all $n$ of the $z_1, \ldots, z_n$
are different.  Then $\Phi * \widehat f$ is compact and
countable, and hence scattered.
\end{lemma}
\begin{proof}
Compactness follows from the compactness of $X,Y$.
By $n$--tightness, $| \widehat f _y | < n$,
and hence $\Phi((  \widehat f _y )^{n}) = \{0\}$,
for all but countably many $y$ (see Lemma \ref{lemma-tight-equiv}).
But for all $y$,  $\pi\iv\{y\}$ is scattered,
so that $\Phi( ( \widehat f _y )^{n})$ is also scattered,
and hence countable.  Thus, the union of all these sets is also countable.
\end{proof}

Dissipation is a notion of smallness, which is balanced
by a notion of bigness, which is really a partition property:

\begin{definition}
\label{def-big}
Fix a real $r > 0$.
The compact space $X$ is $n$--\emph{big} iff for all
$\AA \lea C(X)$ and all partitions $\Upsilon : \AA \to \omega$, there are 
$f_1, \ldots, f_n \in \AA$ and a point $c \in X$
such that
the $\Upsilon(f_j)$, for $j = 1, \ldots, n$, are all equal,
and such that $|f_i(c) - f_j(c)| \ge r$ whenever $1 \le i < j \le n$.
\end{definition}

Since $\AA$ is a linear subspace, it does not matter
which $r > 0$ we use.
The notion of $1$--big is trivial, and
$2$--big is easily characterized:

\begin{lemma}
\label{lemma-two-big}
The compact space $X$ is $2$--big iff $X$ is not second
countable.
\end{lemma}
\begin{proof}
Note that 
$\exists c \,\big(|f_1(c) - f_2(c)| \ge r\big)$ holds iff $\|f_1 - f_2\| \ge r$.
Also, if $X$ is not second countable then $C(X)$ is not separable, and
hence any $\AA \lea C(X)$ is not separable, since
the algebra  generated by the functions in $\AA$
and their complex conjugates is
dense in $C(X)$ by the Stone-Weierstrass Theorem.
\end{proof}

We relate this to the NTIP with the aid of:

\begin{definition}
\label{def-S}
For each $n \ge 2$, define $\Xi_n : \CCC^{n} \to \CCC$ by:
\[
\Xi_n(z_1, \ldots, z_n) = 2 \cdot \prod_{1 \le i < j \le n} (z_j - z_i) \ \ .
\]
\end{definition}

\begin{lemma}
$\Xi_n$ is a polynomial in $n$ variables.
$\Xi(z_1, \ldots, z_n) = 0$ unless all $n$ of the $z_1, \ldots, z_n$
are different. 
If $|z_i - z_j| \ge 1$ for all $i < j \le n$, then either
$|\Re( \Xi_n(z_1, \ldots, z_n)) | \ge 1$ or
$|\Im( \Xi_n(z_1, \ldots, z_n))  | \ge 1$.
\end{lemma}

\begin{lemma}
\label{lemma-big-dis}
Assume that $X$ is compact, $\AA \leac C(X)$, 
$H \subseteq V \subseteq X$, where $H,V$ are both clopen,
and for some $n \ge 2$, 
$V$ is $n$--superdissipated and $H$ is $n$--big.
Assume also that there is a $\psi \in \AA$ such that
$|\psi(x)| \le 1/2$ for all $x \in X \backslash V$ and
$|\psi(x)| \ge 1$ for all $x \in H$.
Then $\AA$ has a non-trivial idempotent.
\end{lemma}

\begin{proof}
Fix $\pi : V  \onto Y$ which is $n$--supertight.
Applying Lemmas \ref{lemma-sing} and \ref{lemma-lim-finer},
we assume also that
we have  $b\in Y$ and $a \in V$ such that $\pi\iv\{b\} = \{a\}$.
Let $\rr^+( z_1, \ldots, z_n ) = -\rr^-( z_1, \ldots, z_n ) =
\Re(\Xi_n( z_1, \ldots, z_n  ))$ and 
$\ii^+( z_1, \ldots, z_n ) = -\ii^-( z_1, \ldots, z_n ) =
\Im(\Xi_n(  z_1, \ldots, z_n ))$, 
so that $\rr^+, \rr^-, \ii^+, \ii^- : \CCC^{n} \to \RRR$.
Call $(E, \rho, \tau)$ \emph{good} iff:
\begin{itemizz}
\item[1.] $\rho, \tau \in \QQQ$ and $1/2 < \rho < \tau < 1$. 
\item[2.] $E \subseteq Y \times \CCC$.  
\item[3.] $[\rho,\tau]$ is disjoint from each of
$\rr^+ * E$, $\rr^- * E$, $\ii^+ * E$, $\ii^- * E$.  
\item[4.] $|\Xi_n( z_1, \ldots, z_n  )| < \rho$ whenever
$z_1, \ldots, z_n\in E_b$. 
\end{itemizz}
For $f\in C(X)$, use $\widehat f$ for $\widehat{\,f \res V\,}$.
Observe that for each $f\in C(X)$, we may choose $\rho,\tau$ so that
$(\widehat f, \rho, \tau)$  is good:
(4) is no problem since $\widehat f _b$ is a singleton.
For the rest, note that each of
$\rr^+ * \widehat f$, $\rr^- * \widehat f$, $\ii^+ * \widehat f$,
$\ii^- * \widehat f$
is scattered by Lemma \ref{lemma-scat-phi},
so we may choose $\rho,\tau$ to make (1)(3) true.

Let $\BB$ be a countable open base for $Y \times \CCC$
which is closed under finite unions.
For each $f \in \AA$, choose $s = s_f \in \omega$ so that
$|(\psi(x))^s \, f(x)| \le 1/8$ for all $x\in X \backslash V$.
Then, choose $\rho_f, \tau_f$ so that
$(\widehat{\psi^{s_f} f\,}, \rho_f, \tau_f)$
is good.  Then, applying Lemma \ref{lemma-open-cover}, 
choose a $U_f \in \BB$ such that
that $(U_f,\rho_f,\tau_f)$ is good and
$\widehat{\psi^{s_f} f\,} \subseteq U_f$.
Next, apply the definition of ``$n$--big'' using $\AA \res H$:
Fix $c \in H $ and
$f_1,f_2, \ldots , f_n \in \AA$ and $(U,\rho,\tau,s)$
such that $(U_{f_j},\rho_{f_j},\tau_{f_j}) =  (U,\rho,\tau)$
and $s_{f_j} = s$ for all $j$,
and also $|f_j(c) - f_k(c)| \ge 1$,
and hence $|\psi^s(c) f_j(c) - \psi^s(c) f_k(c)| \ge 1$,
whenever $j \ne k$.

Let $h(x) = \Xi_n( (\psi(x))^s \, f_1(x), \ldots, (\psi(x))^s \, f_n(x))$;
then $h \in \AA$.
Then, choose $\Phi \in \{\rr^+, \rr^-, \ii^+, \ii^-\}$ so that
$\Phi(  (\psi(c))^s f_1(c), \ldots, (\psi(c))^s f_n(c)) \ge 1$,
and let $k(x) = \Phi( (\psi(x))^s \,f_1(x), \ldots,(\psi(x))^s \,f_n(x))$;
so $k(x)$ is either $\pm\Re(h(x))$ or  $\pm\Im(h(x))$.

Note that when $x \in X \backslash V$, 
each $|(\psi(x))^s \, f_j(x) - (\psi(x))^s \, f_k(x))| \le 1/4$
so (referring to the definition of $\Xi$),
$|h(x)| \le 1/2$.
Then $k(X) = k(V) \cup k(X \backslash V) \subseteq \Phi * U \cup [-1/2,1/2]$
is disjoint from $[\rho,\tau]$, but
contains $k(c) > \tau$ and $k(a) < \rho$.
Thus, either $\Re(h(X))$ or $\Im(h(X))$ is not connected, so 
$\AA$ contains a non-trivial idempotent by Lemma \ref{lemma-get-clopen}.
\end{proof}

In this section, we use only the special case of this lemma where
$H = V = X$, in which case the hypotheses on $\psi$ are trivial,
and the above proof can be simplified somewhat.  The more general
result will be needed in Section \ref{sec-remov}.

Setting $H = V = X$, we have:

\begin{lemma}
\label{lemma-big-ntip}
Suppose that $n \ge 2$ and  $X$ is both $n$--big and $n$--superdissipated.
Then $X$ has the NTIP.
\end{lemma}

Applying this and Lemma \ref{lemma-two-big}, we have:

\begin{theorem}
If $X$ is $2$--superdissipated and
is not second countable, then $X$ has the NTIP.
\end{theorem}

This theorem yields the NTIP for some spaces not covered by
\cite{HK1, KU}, but the result on CSWP, obtained from Corollary \ref{cor-idem},
is contained in the results of \cite{HK1}:

\begin{corollary}
\label{cor-CSWP-onelim}
If $X$ is $2$--superdissipated and
does not contain a Cantor subset, then $X$ has the CSWP.
\end{corollary}

The examples of \cite{HK2, HK3} show (under $\diamondsuit$ or CH)
that this need not hold if $X$ is merely $2$--dissipated.
To extend this corollary to $3$--\emph{super}dissipated spaces, we
need a mechanism (Lemma \ref{lemma-three-big})
for proving that a space is $3$--big.  This notion,
unlike $2$--big (see Lemma \ref{lemma-two-big}), does not seem to 
have a simple equivalent in terms of standard cardinal functions;
see Section \ref{sec-rem}.  

\begin{lemma}
\label{lemma-get-good}
Assume that $n \ge 1$ and that $X$ is
$(n+2)$--superdissipated but not $(n+1)$--superdissipated,
and then fix $\sigma: X \onto Z$ which is $(n+2)$--supertight,
where $Z$ is compact metric.
Assume that $X$ does not have a Cantor subset.
Fix $\AA \lea C(X)$ and $\Upsilon : \AA \to \omega$.
Fix any disjoint open sets $V_0,V_1,V_2 \subseteq \CCC$
and any $\pi \in \MMMM(X)$.
Then there are $f,g,a,d,c$ such that:
\begin{itemizz}
\item[1.] $f,g\in \AA$ and
$\Upsilon(f) = \Upsilon(g)$.
\samepage  
\item[2.] $a, d  \in X$, 
$c = \sigma(a) = \sigma(d) \in Z$, and
$\pi(a) = \pi(d)$.
\item[3.] $f(a) \in V_0$ and $g(a) \in V_1$.
\item[4.] $f(d) \in V_2$ and $g(d) \in V_2$.
\item[5.] For all $x \in \sigma\iv\{c\}$,
$\ (f(x),g(x)) \in V_0 \times V_1 \;\cup\; V_2 \times V_2$.
\end{itemizz}
\end{lemma}
\begin{proof}
First, replacing $\pi$ by a finer map, we may assume that 
$\pi \le \sigma$, so that $\pi \in \MMMMS(X)$ 
and $\pi$ also is $(n+2)$--supertight
(see Lemmas  \ref{lemma-MS-down}, \ref{lemma-count}, and \ref{lemma-lim-finer}).
Say $\pi: X \onto Y$;
then fix $\Gamma \in C(Y,Z)$ with $\sigma = \Gamma \circ \pi$.

Since $\pi$ is not $(n+1)$--supertight, fix a loose family for $\pi$,
$\{P_0, \ldots, P_n\}$, with each $\pi(P_j) = Q$
and each $P_j$ perfect (see Lemma \ref{lemma-tight-equiv}).
Then $\{P_0, \ldots, P_n\}$ is also a loose family for
$\sigma$, with each $\sigma(P_j) = \Gamma(Q)$;
note that $\Gamma(Q)$ cannot be scattered since
$Q$ is not scattered and each $\Gamma\iv\{z\}$ is scattered.
Then $\sigma\iv(\Gamma(Q)) = \pi\iv(\Gamma\iv(\Gamma(Q)))$ is 
superdissipated by Lemma \ref{lemma-get-one},
so it has the CSWP by Corollary \ref{cor-CSWP-onelim}.
Also, $X$ is totally disconnected by Lemma \ref{lemma-discon}.
Fix closed disjoint $\widetilde{P_j} \subseteq \sigma\iv(\Gamma(Q))$
such that each $\widetilde{P_j} \supseteq P_j$ and
$\bigcup_j P_j = \sigma\iv(\Gamma(Q))$.
Note that each $\sigma \res \widetilde{P_j}$ is supertight by
Lemma \ref{lemma-get-one}.

Choose $y_\xi \in Q$ for $\xi < \omega_1$ such that
the $\Gamma(y_\xi)$ are all different and each
$|\pi\iv\{y_\xi\} \cap P_0| \ge 2$; this is possible because
$P_0$ does not have a Cantor subset.
Then, applying the CSWP for  $\sigma\iv(\Gamma(Q))$,
choose $h_\xi \in \AA$ for $\xi < \omega_1$ such that
$h_\xi(\widetilde{P_0}) \subseteq V_0 \cup V_1$,
$h_\xi(\widetilde{P_j}) \subseteq V_2$ when $j \ge 1$,
and $h_\xi( \pi\iv\{y_\xi\} \cap P_0 )$ meets both $V_0$ and $V_1$.
Since there are only countably many values for $\Upsilon$,
we may assume that the $\Upsilon(h_\xi)$ are all the same.
For each $\xi$, we have $P_0$ partitioned into two relatively
clopen sets, $h_\xi\iv(V_0) \cap P_0$ and $h_\xi\iv(V_1) \cap P_0$,
and both these sets meet $\pi\iv\{y_\xi\}$.
If these clopen partition were the same for all $\xi$, we would
contradict the tightness of $\pi\res P_0$ (see Lemma \ref{lemma-get-one}),
so that we may fix $\xi \ne \eta$ with
$H := h_\xi\iv(V_0) \cap h_\eta\iv(V_1) \cap P_0$ non-empty,
and thus perfect. Let $f = h_\xi$ and $g = h_\eta$.

If we choose any $a \in H$, we may set 
$c = \sigma(a)$, and choose any $d \in P_1 \cap \pi\iv\{\pi(a)\}$.
This will satisfy (1)(2)(3)(4), but (5) might fail, since there may
be an $x_c \in \sigma\iv\{c\}$ such that $x_c \in \widetilde{P_0}$
and either $f(x_c) \in V_1$ or  $g(x_c) \in V_0$.
But note that we also have $a \in \sigma\iv\{c\}$ and
$a \in \widetilde{P_0}$ and $f(a) \in V_0$ and  $g(a) \in V_1$.
Consider the map $(f,g) : X \to \CCC\times\CCC$.
If (5) fails for every choice of $a \in H$, then there would be
uncountably many $c \in \pi(H)$ such that $(f,g)$ takes more than
one value on $\widetilde{P_0} \cap \sigma\iv\{c\}$, contradicting 
the tightness of $\sigma\res \widetilde{P_0}$.
Thus, we may choose $a,b,d$ so that (1)(2)(3)(4)(5) hold.
\end{proof}

\begin{lemma}
\label{lemma-three-big}
Assume that $X$ is not dissipated,
but that $X$ is $m$--superdissipated for some $m \in \omega$,
and that $X$ does not have a Cantor subset.
Then $X$ is $3$--big.
\end{lemma}
\begin{proof}
Fix $\AA \lea C(X)$ and $\Upsilon : \AA \to \omega$.
Fix any disjoint open sets $V_0,V_1,V_2 \subseteq \CCC$.
To verify that $X$ is $3$--big, it is sufficient to find
$h_0,h_1,h_2 \in \AA$ and $x \in X$ such that each $h_j(x) \in V_j$.

Fix $n \ge 1$  such that $X$ is
$(n+2)$--superdissipated but not $(n+1)$--super\-dissipated,
and then fix $\sigma: X \onto Z$ which is $(n+2)$--supertight
Let $\BB$ be a countable open base for $Z$.
For $\pi \in \MMMM(X)$, call 
$F = (f, g, a, d, c ,s, U) = (f^F, g^F, a^F, d^F, c^F ,s^F, U^F)$
\emph{good for} $\pi$ iff (1--5) from Lemma \ref{lemma-get-good} 
hold together with:

\begin{itemizz}
\item[6.] $s \in \omega$ and $\Upsilon(f) = \Upsilon(g) = s$.
\item[7.] $c \in U$,  $U \in \BB$,  and for all $x \in \sigma\iv(\overline U)$,
$\ (f(x),g(x)) \in V_0 \times V_1 \;\cup\; V_2 \times V_2$.
\end{itemizz}
Such an $F$ always exists.  To see this, first get
$(f, g, a, d, c)$ by Lemma \ref{lemma-get-good} to satisfy (1--5).
Then (6) is trivial, and we choose $U$ to satisfy (7) using the fact that
$\{z \in Z : \forall x \in \sigma\iv\{z\}\,
[(f(x),g(x)) \in V_0 \times V_1 \;\cup\; V_2 \times V_2]\}$
is open.

Note that if $F$ is good for $\pi$ and $\pi \le \varphi$
then $F$ is good for $\varphi$.

Next, note that there are fixed $s$ and $U$ such that
for all $\pi \in \MMMM(X)$, there is an $F$ good for $\pi$ with
$s^F = s$ and $U^F = U$:
If not, then for each $s,U$, choose
$\varphi_{s,U}$ such that no $F$ good for $\varphi$ satisfies $s^F = s$
and $U^F = U$.
Then fix $\pi$ such that $\pi \le \varphi_{s,U}$ for each $s,U$.
An $F$ which is good for $\pi$ yields a contradiction.

For each $\pi$, choose $F^\pi$ good for $\pi$ 
with $s^{F^\pi} = s$ and $U^{F^\pi} = U$,
and write
$(f^\pi, g^\pi, a^\pi, d^\pi, c^\pi)$ for
$(f^{F^\pi}, g^{F^\pi}, a^{F^\pi}, d^{F^\pi}, c^{F^\pi})$.

Now, for each $\pi$, we have $\sigma\iv(\overline U)$
partitioned into two relatively
clopen sets, $A^\pi = \{x \in \sigma\iv(\overline U) :
(f^\pi(x), g^\pi(x)) \in V_0 \times V_1\}$ and
$D^\pi =
\{x \in \sigma\iv(\overline U) :
(f^\pi(x), g^\pi(x)) \in V_2 \times V_2\}$.
If these are all the same, say
$A^\pi = A$ and $D^\pi = D$ for all $\pi$;
then we may fix $\pi\in C(X,[0,1])$ which is $0$ on $A$ and $1$ on $D$,
so $\pi(a^\pi) = 0$ and $\pi(d^\pi) = 1$, contradicting (2).
Thus, we can choose $\pi, \varphi$
and an $x \in A^\pi \cap D^\varphi$; then
$f^\pi(x)  \in V_0$, $g^\pi(x)  \in V_1$,
$f^\varphi(x)  \in V_2$, as required.
\end{proof}

The ``obvious'' generalization of this would say that if 
$X$ does not have a Cantor subset and is
$(n+2)$--superdissipated but not $(n+1)$--superdissipated,
then $X$ is $(n+2)$--big.
For $n=1$ this is Lemma \ref{lemma-three-big}, and
for $n=0$ this is Lemma \ref{lemma-two-big}.
Unfortunately, this is not true in general; see 
Example \ref{ex-big-limprop}.  We do get:

\begin{theorem}
\label{thm-big-cswp}
Assume that $X$ is compact and is $3$--superdissipated and
does not have a Cantor subset.  Then $X$ has the CSWP.
\end{theorem}
\begin{proof}
Since ``$3$--superdissipated'' is closed-hereditary, it is sufficient, by 
Corollary \ref{cor-idem}, to 
assume that $X$ is also perfect and prove that $X$ has the NTIP.
$X$ cannot be second countable, so $X$ is 
$2$--big by Lemma \ref{lemma-two-big}.
If $X$ is not $2$--superdissipated, then $X$ is
$3$--big by Lemma \ref{lemma-three-big}.
Thus, whether or not  $X$ is $2$--superdissipated,
it has the NTIP by Lemma \ref{lemma-big-ntip}.
\end{proof}

\begin{corollary}
\label{cor-two-lots}
If $X$ is compact and $X \subseteq L_0 \times L_1$,
where $L_0,L_1$ are a LOTSes, then $X$ has the CSWP
iff $X$ does not contain a copy of the Cantor set.
\end{corollary}
\begin{proof}
By Lemma \ref{lemma-one-S}, we may assume that $X \subseteq (I_S)^2$.
Then $X$ is $3$--superdissipated by Lemma \ref{lemma-big-lim},
so $X$ has the CSWP by Theorem \ref{thm-big-cswp}.
\end{proof}

We now can extend this to products of three LOTSes, using
an argument which is much more specific to ordered spaces.
First, we introduce a notation for lines, boxes, etc.\@ in such products.

\begin{definition}
\label{def-boxes}
Let $\prod_{\alpha < \kappa} L_\alpha$ be a product of LOTSes,
and use $<$ for the order on each $L_\alpha$.  Then:
\begin{itemizn}{"2B}
\item
If $\beta < \kappa$ and $c$ is a point in
$\prod_{\alpha \ne \beta} L_\alpha$,
then $\lin(\beta, c) = \{x \in \prod_{\alpha < \kappa} L_\alpha :
\forall \alpha \ne \beta \, [ x_\alpha = c_\alpha ]\}$. 
A \emph{line} in  $\prod_{\alpha < \kappa} L_\alpha$
is any set of the form $\lin(\beta, c)$.
\item
$<^+$ is $<$\textup; \ $<^-$ is $>$\textup;
\ $\le^+$ is $\le$\textup;\  $\le^-$ is $\ge$ .
\item
$\DD = \{+, -\}^\kappa$ is the set of all \emph{directions}.
For $\Delta \in \DD$ and $x,y \in \prod_{\alpha < \kappa} L_\alpha$,
$x <^\Delta y$ iff
$\forall \alpha \, [ x_\alpha <^{\Delta_\alpha}y_\alpha]$ and
$x \le^\Delta y$ iff
$\forall \alpha \, [ x_\alpha \le^{\Delta_\alpha}y_\alpha]$.
\item 
If $a,b \in \prod_{\alpha < \kappa} L_\alpha$, 
then
$\bx[a,b] = \prod_{\alpha < \kappa}
[\min(a_\alpha,b_\alpha), \max(a_\alpha,b_\alpha)] $,
and a \textup(closed\textup) \emph{box} is any set of this form.
\item 
If $a \in \prod_{\alpha < \kappa} L_\alpha$ and $\Delta \in \DD$, then
$\corn(a,\Delta) = \{x \in \prod_{\alpha < \kappa} L_\alpha 
: a \le^\Delta x\}$.
\item 
If $a \in B \subseteq \prod_{\alpha < \kappa} L_\alpha$ and
$\Delta \in \DD$, then
$\corn(a, B, \Delta) = B \cap \corn(a,\Delta)$.
\end{itemizn}
\end{definition}
For example, in $\RRR^3$:
$(2,4,6) <^{+-+} (3,3,7) \le^{+-+} (4,2,7)$.  Now,
let
$B =  [0,9]^3 = \bx[(0,0,0), (9,9,9)] = \bx[(9,0,9), (0,9,0)]$.  Then
$\corn( (2,4,6), B, \mop\mom\mop)$ is the box
$[2,9]\times [0,4]\times [6,9]$.
The directions $\Delta\in\DD$ are also useful inside products
of the form $(I_S)^\kappa$.  Continuing the notation of
Definition \ref{def-das},

\begin{definition}
If $\sigma: (I_S)^\kappa \onto I^\kappa$ is the standard map,
$y \in I^\kappa$, and $\Delta\in \DD$, then
$y^\Delta = \langle y^{\Delta_\alpha} : \alpha < \kappa\rangle$.
\end{definition}
For example, 
if $b = (b_0,b_1,b_2) \in I^3$, then $\sigma\iv\{b\}$ consists of 
the points, $b^{\pm\pm\pm} = (b_0^\pm,b_1^\pm,b_2^\pm)$;
e.g., $b^{+-+}$ denotes the point
$(b_0^+,b_1^-,b_2^+) \in (I_S)^3$.
The size of $\sigma\iv\{b\}$ will be $8,4,2$ or $1$ depending on whether
$3,2,1$ or $0$ of the $b_0,b_1,b_2$ lie in $S$.

The following lets us establish bigness
for subsets of $(I_S)^n$ by checking a simpler
geometric property:

\begin{lemma}
\label{lemma-big-helper}
Fix $S \subseteq (0,1)$, a closed
$X \subseteq (I_S)^n$, and $m$ with $2^{n-1} < m \le 2^n$.
Assume that whenever $\Upsilon : S^n \to \omega$, there are distinct
$ \Delta_1,  \Delta_2, \ldots,  \Delta_m \in \DD$, a point $x \in X$,
and
$ d_1,  d_2, \ldots,  d_m \in S^n$ such that
$x \in \corn(d_j^{\Delta_j}, \Delta_j)$ for each $j$,
and such that $ \Upsilon(d_1) =  \Upsilon(d_2) = \cdots =  \Upsilon(d_m)$.
Then $X$ is $m$--big.
\end{lemma}
\begin{proof}
Note that for $d  \in S^n$, the points $d^\Delta \in (I_S)^n$,
for $\Delta \in \DD$, are all distinct,
and the $\corn(d^\Delta, \Delta)$, for $\Delta \in \DD$,
partition $(I_S)^n$ into $2^n$ clopen subsets.

Fix $\AA \lea C(X)$ and $\Upsilon : \AA \to \omega$.
Since finite spaces have the CSWP, we may choose,
for each $d \in S^n$, an $f_d \in C((I_S)^n)$ with $f_d \res X \in \AA$
such that
the $f_d(d^\Delta)$, for $\Delta \in \DD$,
are $2^n$ distinct integers.
We shall verify the definition of ``$m$--big'' just by considering
the functions $f_d \res X$; the $r$ in Definition \ref{def-big}
will be $1/2$.

Each $f_d$ is continuous, so choose $p(d), q(d) \in \QQQ^n$
with $\forall \mu\,[ p(d)_\mu < d_\mu < q(d)_\mu]$  such that
$\sup\{ | f_d(x) - f_d(d^\Delta) | :
x \in \corn(d^\Delta,\; \bx[p(d)^{+},q(d)^{-}]  ,\; \Delta)\} \le 1/4$
for each $\Delta \in \DD$.  Here, for $y \in I^n$,
$y^+$ abbreviates $(y_0^+, \ldots, y_{n-1}^+)$ and
$y^-$ abbreviates $(y_0^-, \ldots, y_{n-1}^-)$.
Now, let $\Upsilon'(d) = (\Upsilon(f_d \res X),\; p(d),\; q(d))$.
Since $\ran(\Upsilon')$ is countable, we may apply the hypotheses of
the lemma and fix distinct
$ \Delta_1,  \Delta_2, \ldots,  \Delta_m \in \DD$, along with $x \in X$
and $ d_1,  d_2, \ldots,  d_m \in S^n$, such that
$x \in \corn(d_j^{\Delta_j}, \Delta_j)$ for each $j$,
$\; \Upsilon(f_{d_1}\res X) =
 \Upsilon(f_{d_2}\res X) = \cdots =  \Upsilon(f_{d_m}\res X)$,
and also each $p(d_j) = p$ and $q(d_j) = q$ for some $p,q \in \QQQ^n$.

Note that all $d_j^\Delta \in \bx[p^{+},q^{-}]$.  Also, since $m > 2^{n-1}$,
$\{ \Delta_1,  \Delta_2, \ldots,  \Delta_m\}$ contains both
$\Delta$ and $-\Delta$ for some $\Delta$, which implies
(using $x \in \corn(d_j^{\Delta_j}, \Delta_j)$) that
$x \in \bx[p^{+},q^{-}]$.  Thus,
$x \in \corn(d_j^{\Delta_j},\; \bx[p^{+},q^{-}],\; \Delta_j)$, so
$| f_{d_j}(x) - f_{d_j}(d_j^{\Delta_j}) | \le 1/4$ for each $j$,
so that  $|f_{d_j}(x) - f_{d_k}(x)| \ge 1/2$ when $j \ne k$.
\end{proof}

Note that the points $ d_j^{\Delta_j}$ were not assumed to lie in $X$.

\begin{lemma}
\label{lemma-six-big}
Assume that $S \subseteq (0,1)$,
$X$ is a closed subspace of $(I_S)^3$, $X$ is not $3$--dissipated,
and $X$ does not contain a Cantor subset.
Then $X$ is $6$--big.
\end{lemma}
\begin{proof}
We verify the hypotheses of Lemma \ref{lemma-big-helper},
so fix $\Upsilon : S^3 \to \omega$; we must find appropriate
$ \Delta_1,  \Delta_2, \ldots,  \Delta_6 \in \DD = \{+, -\}^3$,
$x \in X$, and $ d_1,  d_2, \ldots,  d_6 \in S^3$.

Note that it is sufficient to find $x$ along with 
points $c_E, c_F, c_G \in S \times S$, numbers
$u_E, w_E, u_G, w_G, u_F, v_F, w_F \in S$, and $\Delta \in \{+, -\}^2$
such that:
\begin{itemizz}
\item[1.] $c_E <^\Delta c_F <^\Delta c_G$.
\item[2.] $u_E, u_F, u_G < v_F$ and $v_F < w_E, w_F, w_G$.
\item[3.] $\Upsilon$ has the same value on the $6$ points:
$
d_1 = (c_E, u_E),
d_2 = (c_E, w_E),\\
d_3 = (c_F, u_F),
d_4 = (c_F, w_F),
d_5 = (c_G, u_G),
d_6 = (c_G, w_G) 
$.
\item[4.] $x \in X$ and $x$ is one of the four points
$(c_F^{ \Gamma}, v_F^\pm)$, where $\Gamma \in \{+, -\}^2$
and $\Gamma$ is different from $\Delta$ and $-\Delta$.
\end{itemizz}
Note that no ordering is assumed among 
$u_E, u_F, u_G$ or among $w_E, w_F, w_G$.
To verify that (1--4) are sufficient, and to clarify our notation,
assume WLOG that $\Delta = ++$, so $c_E <^{++} c_F <^{++} c_G$.
Then $\Gamma$ is either $+-$ or $-+$; WLOG $\Gamma = +-$,
so we are assuming $X$ contains at least
one of the two points $(c_F^{+-}, v_F^\pm)$, denoted by $x$. 
But now we obtain the hypotheses of
Lemma \ref{lemma-big-helper}.
Namely,  $x \in \corn(d_j^{\Delta_j},  \Delta_j)$ for 
$j = 1,2,\ldots,6$, setting
\  $
 \Delta_1 = \mop\mop\mop \,,
 \, \Delta_2 = \mop\mop\mom \,,
 \, \Delta_3 = \mop\mom\mop \,,
 \, \Delta_4 = \mop\mom\mom \,, \allowbreak
 \, \Delta_5 = \mom\mom\mop \,,
 \, \Delta_6 = \mom\mom\mom 
$

Now, to obtain (1--4):
If $E \subseteq S$, let $\sigma_E: (I_S)^3 \onto (I_E)^3$ be the natural
map; so $\sigma_\emptyset = \sigma$.
If also $E \in [S]^\omega$ (i.e., $|E| = \aleph_0$),
then $(I_E)^3$ is a compact metric space,
and we shall use the fact that \emph{none} of these
$\sigma_E$ are $3$--tight.

If $E_1 \subseteq E_2 \in [S]^\omega$ then
$\sigma_{E_2} \le \sigma_{E_1}$ (see Lemma \ref{lemma-fine-equiv}).
Observe that $ [S]^\omega $  is countably directed upward.
Call $\UU \subseteq [S]^\omega$ \emph{cofinal} iff
$\forall E_1 \in [S]^\omega \, \exists E_2 \in \UU
\,  ( E_1 \subseteq\nobreak E_2 )$;
then $\UU$ is also countably directed upward.
We shall use this observation
several times to show that a number of quantities dependent on $E$
can in fact be chosen uniformly, independently of $E$, on a cofinal set.

Temporarily fix an $E  \in [S]^\omega$.
Then we have $P_j = P^E_j \subseteq X \subseteq (I_S)^3 $
for $j = 0,1,2$ such that $\{P_0, P_1, P_2\}$ is a loose family.
Then each $\sigma_E(P_j) = Q$,
where $Q  = Q^E\subseteq \sigma_E(X) \subseteq (I_E)^3$ is uncountable.
We can now get such a $Q$ to be of a very simple form:

First, note that $Q$ must be a subset of finite union of lines.
If not, then we may choose $y^\ell = (y^\ell_0, y^\ell_1, y^\ell_2) \in Q$
for $\ell \in \omega$ such that no two of the $y^\ell$ lie
on the same line; that is, whenever $\ell < m < \omega$,
the triples  $y^\ell$  and $y^m$ differ on at least two coordinates.
Now, we may thin the sequence and permute the coordinates and assume
that each of the two sequences
$\langle y^\ell_0 : \ell \in \omega \rangle$
and
$\langle y^\ell_1 : \ell \in \omega \rangle$
is either strictly increasing or strictly decreasing, while
$\langle y^\ell_2 : \ell \in \omega \rangle$ is either constant or
strictly increasing or strictly decreasing.
If $H$ is the set of limit points of the sequence
of sets $\langle  \sigma\iv\{y^\ell_i\} : \ell \in \omega \rangle$,
then $|H| \le 2$, but $H$ must meet each of $P_0,P_1,P_2$,
which is a contradiction.

Next, shrinking $Q$, along with $P_0,P_1,P_2$,
we may assume that $Q = Q^E$ is a subset of one line;
say $Q^E \subseteq \lin(\beta_E, c_E)$, where $\beta_E < 3$.

$\beta_E$ depends on $E$, but since $ [S]^\omega $  is
countably directed upward, there is a fixed $\beta$ such
that $\beta_E = \beta$ on a cofinal set 
$\UU \subseteq [S]^\omega$.  By permuting coordinates, we may
assume $\beta = 2$, so that  $Q^E \subseteq \lin(2, c_E) \subseteq (I_E)^3$,
where $c_E = (a_E, b_E) \in (I_E)^2$.
From now on, we shall delete the ``2''; so
$\lin(c_E) = \{(a_E, b_E,u) : u \in I_E \}$.
Then $Q^E = \{c_E\} \times \widetilde Q_E $, where
$\widetilde Q_E \subseteq I_E$.

Again, fix $E$, and temporarily delete some of the sub/super-script $E$.
Now $\sigma_E\iv(\lin( c)) \subseteq (I_S)^3$ is a union of 
1, 2, or 4 lines in $(I_S)^3$.  However, the existence
of $Q,P_0,P_1,P_2$ implies that
$\sigma_E : \sigma_E\iv(\lin(c)) \onto \lin(c)$ is not $3$--tight,
so in fact $\sigma_E\iv(\lin( c))$
is a union of 4 lines, which means that $a,b \in S \backslash E$;
that is, we may regard $a,b$ as real numbers which are not split
in $I_E$, but which are split into $a^\pm, b^\pm$ in $I_S$,
and $\sigma_E\iv(Q) \subseteq
\lin( c^{++} )\cup 
\lin( c^{+-} )\cup 
\lin( c^{-+} )\cup 
\lin( c^{--} )
\subseteq (I_S)^3$.
Now $\sigma_E\iv(Q) \cap \lin( c^{++} ) \cap X$ is some closed subset of
$\sigma_E\iv(Q) \cap \lin( c^{++} )$, but replacing $Q$ by a
smaller perfect set, we may assume that this closed subset is either
empty or all of $\sigma_E\iv(Q) \cap \lin( c^{++} )$.
Repeating this argument three more times, we may assume that 
each of the four sets $\sigma_E\iv(Q) \cap \lin( c^{\pm \pm } )$
is either contained in $X$ or disjoint from $X$.
Again, the existence of $P_0,P_1,P_2$ implies that
$\sigma_E : \sigma_E\iv(Q) \cap X \onto Q$ is not $3$--tight,
so at least three of the four sets $\sigma_E\iv(Q) \cap \lin( c^{\pm \pm } )$
are contained in $X$.
Which three or four depends on $E$; there is a cofinal
set on which it is the same, although this is irrelevant now.
More importantly, since
$\widetilde Q_E \subseteq I_E$ and $E$ is countable, we may
shrink $Q^E$ and assume that $\widetilde Q_E \cap E = \emptyset$;
that is, we may regard $\widetilde Q_E$ as a perfect subset
of $I \backslash E$.
Note that $S$ must meet every perfect subset of $\widetilde Q$,
since otherwise $X$ would contain a Cantor subset.
In particular, $S \cap \widetilde Q$ is uncountable.
Now $c_E = c = (a,b)$ is fixed, and
for each $u \in S \cap \widetilde Q$, we have the triple 
$d = d_u = (a,b,u)$.
We may now choose  $t\in \omega$ and $u,v,w \in  S \cap \widetilde Q$
such that $u < v < w$ and
$\Upsilon(d_u) = \Upsilon(d_v) = \Upsilon(d_w)  = t$.
Also choose rational $\rho,\tau$ with $u < \rho < v < \tau < w$.

Of course, $t,\rho,\tau,u,v,w$ depend on $E$, but there are only
$\aleph_0$ possibilities for $t,\rho,\tau$, so 
we may assume that for $E$ in our cofinal
set $\UU$, these are always the same, whereas $u,v,w$ are really $u_E,v_E,w_E$.

Choose an increasing $\omega_1$ sequence
$\langle E_\xi : \xi < \omega_1 \rangle$
of elements of $\UU$ such that 
$\xi < \eta \to c_{E_\xi} \in (E_\eta)^2$.
Now $c_{E_\xi} = ( a_{E_\xi} b_{E_\xi})$ and
$ a_{E_\xi},  b_{E_\xi} \notin E_\xi $,
so $a_{E_\eta} \ne a_{E_\xi}$ and
$b_{E_\eta} \ne b_{E_\xi}$ whenever $\xi \ne \eta$.
It follows that we may find distinct 
$\xi_n < \omega_1$ for $n \in \omega$ and a fixed $\Delta \in \{+, -\}^2$
such that $m < n \to c_{E_{\xi_m}} <^\Delta c_{E_{\xi_n}}$.
But, we only need three of these, so let $E,F,G$ denote
$E_{\xi_0}, E_{\xi_1}, E_{\xi_2}$.  Then we have
$c_E <^\Delta c_F <^\Delta c_G$ as in (1) above.
$u_E, u_F, u_G < \rho < v_F < \tau < w_E, w_F, w_G$, so (2) holds.
(3) holds because $\Upsilon$ has the same value $t$ on all
$( a_{E_\xi},  b_{E_\xi} , u_{E_\xi})$,
$( a_{E_\xi},  b_{E_\xi} , v_{E_\xi})$,
$( a_{E_\xi},  b_{E_\xi} , w_{E_\xi})$.
Finally,  we may choose $x$ to make (4) hold because
at least three of the four sets
$\sigma_F \iv(\widetilde Q _ F) \cap \lin( c_F^{\Gamma} )$
(for $\Gamma \in \{+, -\}^2$) 
are contained in $X$ and $v_F \in S \cap \widetilde Q _ F $,
and for these $\Gamma$, both points $(c_F^{ \Gamma}, v_F^\pm)$
lie in $X$.
\end{proof}

\begin{proofof}{Theorem \ref {thm-three-lots}}
By Lemma \ref{lemma-one-S}, we may assume that $X \subseteq (I_S)^3$.
Since the properties assumed of $X$ are
closed-hereditary, it is sufficient, by 
Corollary \ref{cor-idem}, to 
assume that $X$ is also perfect and prove that $X$ has the NTIP.
Note that ``dissipated'' is the same as ``superdissipated'' for
these spaces.
If $X$ is $3$--dissipated, then $X$ has the CSWP, and
hence the NTIP, by Theorem \ref{thm-big-cswp}.
If $X$ is not $3$--dissipated, then $X$ 
is $5$--big by Lemma \ref{lemma-six-big}, but it also 
is $5$--dissipated by Lemma \ref{lemma-big-lim},
so $X$ has the NTIP by Lemma \ref{lemma-big-ntip}.
\end{proofof}

We do not know if the same theorem holds when $X$ is contained in
a product of four LOTSes, but the analogue of Lemma \ref{lemma-six-big}
is false.  That is, there is (see Example \ref{ex-seven})
a closed $X \subseteq (I_S)^4$
such that $X$ is not $8$--dissipated and is not
$7$--big.  Of course, $X$ must be $9$--dissipated,
but to prove the NTIP by our methods, $X$ would need to be
$9$--big.

\section{Removable Spaces}
\label{sec-remov}
The property of a compact space being \emph{removable}, defined in \cite{HK1},
is a strengthening of the CSWP.  Many of the spaces
proved in Section \ref{sec-cswp} to have the CSWP are in fact removable.
We recall the definition, which is in terms of the \v Silov boundary:

\begin{definition}
\label{def-silov}
If $\AA\leac C(X)$, then $\Sh(\AA)$ denotes the \emph{\v Silov boundary};
this is the smallest non-empty closed $H \subseteq X$ such that 
$\|f\| = \sup\{|f(x)| : x \in H\}$ for all $f\in \AA$.
\end{definition}

This is discussed in texts on function algebras; see
\cite{GAM, GA}.
Note that $\Sh(\AA)$ cannot be finite
unless $X$ is finite, in which case $\Sh(\AA) = X$.

\begin{definition}
\label{def-remove}
A compact space $K$ is \emph{removable} iff for all $X,U,\AA$, if:
\begin{itemizn}{"2B}
\item $X$ is compact, $U \subsetneqq X$, and $U$ is open,
\item $\overline U$ is homeomorphic to a subspace of $K$, and
\item $\AA \leac C(X)$ and all idempotents of $\AA$ are trivial,
\end{itemizn}
then $\Sh(\AA) \subseteq X\backslash U$.
\end{definition}

The next four lemmas are clear from \cite{HK1}:

\begin{lemma}
If $X$ is removable, then $X$ is totally disconnected and has the CSWP.
\end{lemma}

It is unknown whether the converse to this lemma is true.
The removable spaces are of interest because
one can prove some theorems about them
which are currently unknown for the CSWP spaces.  In particular,
the removable spaces form a local class
(see Definition \ref{def-local}); this follows from:

\begin{lemma}
If the compact $X$ is a finite union of closed sets, each of which
is removable, then $X$ is removable.
\end{lemma}

More generally,
one can do a type of Cantor-Bendixson analysis for a compact $X$,
iteratively deleting open sets with removable closures; if one
gets to $\emptyset$, then $X$ itself is removable and hence
has the CSWP (see \cite{HK1}, Lemma 2.15).
This results in the next definition and lemma.

\begin{definition}
A compact space $P$ is \emph{nowhere removable} iff
$\overline W$ is not  removable
for all non-empty open $W \subseteq P$.
\end{definition}

\begin{lemma}
If $X$ is compact and not removable, then there is a non-empty
closed $P \subseteq X$ such that $P$ is nowhere removable.
\end{lemma}

In particular, since the one-point space is removable,

\begin{lemma}
Every compact scattered space is removable.
\end{lemma}

\begin{definition}
\label{def-R}
$\RRRR$ is the class of all compact spaces $X$ such
that for all perfect $H \subseteq X$:
There is non-empty relatively clopen $U \subseteq H$ such that 
either $U$ is removable or
for some finite $n \ge 2$, $U$ is both $n$--big and $n$--superdissipated.
\end{definition}

If $X$ is removable, then $X \in \RRRR$, and we shall soon prove the 
converse statement.
No space in $\RRRR$ can contain a Cantor subset
(since the Cantor set is neither $2$--big nor removable).
All spaces in $\RRRR$ are totally disconnected by Lemma \ref{lemma-discon}.

Our proof will use the following restatement of Definition \ref{def-remove}:

\begin{lemma}
\label{lemma-restate}
Assume that $\KKKK$ is a closed-hereditary class of
totally disconnected compact spaces,
and assume that whenever $Z,V,\AA$ satisfy:
\begin{itemizn}{"2B}
\item $Z$ is compact and infinite, $\AA \leac C(Z)$, and $\Sh(\AA) = Z$.
\item $V \subseteq Z$, $V$ is clopen and non-empty, and $V \in \KKKK$.
\end{itemizn}
then $\AA$ contains a non-trivial idempotent.
Then, all spaces in $\KKKK$ are removable.
\end{lemma}
\begin{proof}
Fix $K \in \KKKK$.  Then fix $X,U,\AA$ satisfying the hypotheses
of Definition \ref{def-remove}.  Let $Z = \Sh(\AA)$.
Assume that
$Z \not\subseteq X \backslash U$.  We shall derive a contradiction.
Shrinking $U$, we may assume that $U$ is clopen.
Clearly $U \ne \emptyset$, so $|X| \ge 2$ (by
$U \subsetneqq X$), so $X$ is infinite (by all idempotents trivial),
so $Z$ is infinite.

$\AA\res Z \leac C(Z)$ and
$\Sh(\AA\res Z) = Z$.  Let $V = Z \cap U$; then $V \ne \emptyset$.
$V \in \KKKK$ because $\KKKK$ is closed-hereditary.
So, $\AA\res Z $ contains a non-trivial idempotent, $f \res Z$,
where $f \in \AA$.  But then $f^2 - f$ is $0$ on $Z$ and hence on $X$,
so $f$ is an idempotent, contradicting 
the hypotheses of Definition \ref{def-remove}.
\end{proof}

\begin{theorem}
\label{thm-R-remov}
$\RRRR$ is the class of all removable spaces.
\end{theorem}

\begin{proof}
Since $\RRRR$ is clearly closed-hereditary, we may apply
Lemma \ref{lemma-restate} to prove that all spaces in $\RRRR$
are removable.  Thus, assume
that $X$ is compact and infinite, $\AA \leac C(X)$, and $\Sh(\AA) = X$,
and $V \subseteq X$ is clopen and non-empty, and $V \in \RRRR$.
We must show that $\AA$ contains a non-trivial idempotent.
We may assume that $V$ is nowhere removable, and in particular perfect,
since otherwise the result is clear from the definition of ``removable''.
Applying the definition of $\RRRR$, whenever $U$ is a non-empty
clopen subset of $V$, there is an $n_U \ge 2$ and a non-empty
clopen $H$ with $H \subseteq U$ and $H$
both $n_U$--big and $n_U$--superdissipated.
Taking a minimal $n_U$ and shrinking $V$,
we may assume that $V$ itself is
$n$--superdissipated, where $n \ge 2$,
and that whenever $U$ is a non-empty
clopen subset of $V$, there is  a non-empty
clopen $H$ with $H \subseteq U$ and $H$
$n$--big.

Since  $X \backslash V$ is not a boundary,
we may fix $\psi \in \AA$ such that $\|\psi\| > 1$ but
$|\psi(x)| \le 1/2$ for all $x \notin V$.
Then fix a non-empty clopen $H \subseteq V$ such that $|\psi(x)| \ge 1$ for all
$x \in H$.  Shrinking $H$,
we may assume that $H$ is $n$--big.
We now get a non-trivial idempotent by Lemma \ref{lemma-big-dis}.
\end{proof}

\begin{corollary}
\label{cor-3remov}
If $X \subseteq (I_S)^3$ is closed and does not
contain a copy of the Cantor set, then $X$ is removable.
\end{corollary}
\begin{proof}
$X \in \RRRR$  by Lemmas 
\ref{lemma-big-lim},
\ref{lemma-six-big}, 
\ref{lemma-three-big}, and
\ref{lemma-two-big}.
\end{proof}

\section{Powers of the Double Arrow Space}
\label{sec-powers}
Here we show that arbitrary finite powers of the double arrow space
$I_{(0,1)}$ are removable, and hence have the CSWP.
This argument works because there is a certain uniformity in the
standard map from $(I_{(0,1)})^k$ onto $I^k$, which is
captured by the next definition:

\begin{definition}
\label{def-superdup}
For $n \ge 1$, $\pi : X \onto Y$ is $n$--\sdup tight iff
for $y \in Y$ and $0 \le j < n$, there are
$K^j_y \subseteq X$ and
$U^j_y \subseteq Y$ satisfying:
\begin{itemizz}
\item[1.] $X,Y$ are compact, $Y$ is metric, and
the Cantor set does not embed into $X$.
\item[2.] For each $y$: The $K^j_y$, for $j < n$, form a clopen
partition of $X$, and each $|K^j_y \cap \pi\iv\{y\}| \le 1$.
\item[3.] For each $j$: $\{(y,z) :  z \in U^j_y\}$ is open in $Y^2$.
\item[4.] For each $y,j$: $\pi\iv(U^j_y) \subseteq K^j_y$.
\item[5.] For each $y,j$: $ K^j_y \setminus \pi\iv(U^j_y)$ is removable.
\end{itemizz}
$X$ is $n$--\sdup dissipated iff $\pi: X \onto Y$ is
$n$--\sdup tight for some $\pi$ and $Y$.
\end{definition}

Some of the $K^j_y$ and $U^j_y$ may be empty,
so ``$n$--\sdup dissipated'' get weaker as $n$ gets bigger.
Note that (2) implies that $|\pi\iv\{y\}| \le n$ for each $y$,
so that $\pi$ is $n$--supertight by Lemma \ref{lemma-lim-prop},
and $X$ is totally disconnected by Lemma \ref{lemma-discon}.
$X$ is $1$--\sdup dissipated iff $X$ is compact and countable.
The class of $n$--\sdup dissipated spaces is closed-hereditary,
since if we have (1 -- 5) and
$\widetilde X$ is a closed subset of $X$, then we also have
(1 -- 5) for $\widetilde X$, using
$\pi\res \widetilde X : \widetilde X \onto \widetilde Y = \pi(\widetilde X)$, 
$\widetilde K^j_y = K^j_y \cap \widetilde X$, and
$\widetilde U^j_y = U^j_y \cap \widetilde Y$.

\begin{lemma}
\label{lemma-das-superdup}
If  $(I_{(0,1)})^{k-1}$ is removable, then
the standard map $\pi : (I_{(0,1)})^k \onto I^k$ is 
$2^k$--\sdup tight.
\end{lemma}
\begin{proof}
As in Definition \ref {def-boxes}, let $\DD = \{+, -\}^k$.
For $y \in I^k$ and $\Delta \in \DD$, 
let $U^\Delta_y = \{z \in I^k : y <^\Delta z\}$, and
let $K^\Delta_y = \{t \in (I_{(0,1)})^{k} : y^\Delta \le^\Delta t\}$.
Then properties (1 -- 4) are easily verified, and (5) holds 
because $ K^j_y \setminus \pi\iv(U^j_y)$ is covered by finitely many
homeomorphic copies of  $(I_{(0,1)})^{k-1}$.
\end{proof}

We shall eventually prove:

\begin{theorem}
\label{thm-superdup-remov}
If $n < \omega$ and $X$ is $n$--\sdup dissipated,
then $X$ is removable.
\end{theorem}

It follows that $X$ is $n$--\sdup dissipated
iff $X$ is removable and there is a $\pi : X \onto Y$ such
that $Y$ is compact metric and each $|\pi\iv\{y\}|  \le n$.
To prove the $\leftarrow$ direction:
In Definition \ref{def-superdup}, take
all $U^j_y = \emptyset$; the $K^j_y$ may simply 
be chosen arbitrarily to satisfy condition (2).
Thus, the notion of ``$n$--\sdup dissipated'' becomes of little
interest, but it was chosen to make the following proof work:

\begin{proofof}{Theorem \ref{thm-das-power}}
Each $(I_{(0,1)})^k$ is in fact removable.
This follows by induction on $k$,
using Lemma \ref{lemma-das-superdup} and
Theorem \ref{thm-superdup-remov}.
\end{proofof}

We shall now prove Theorem \ref{thm-superdup-remov} by showing
that $X \in \RRRR$ (see Definition \ref{def-R}).

\begin{definition}
A compact space $P$ is \emph{nowhere $n$--\sdup dissipated} iff
$\overline W$ is not  $n$--\sdup dissipated
for all non-empty open $W \subseteq P$.
\end{definition}

\begin{lemma}
\label{lemma-get-nowhere}
If $X$ is perfect and $n$--\sdup dissipated, then there is
a non-empty clopen $V \subseteq X$ and an $m$ with $2 \le m \le n$
such that $V$ is $m$--\sdup dissipated and nowhere
$(m-1)$--\sdup dissipated.
\end{lemma}

Theorem \ref{thm-superdup-remov} will follow easily from the next two
lemmas, about spaces which are $n$--\sdup dissipated and nowhere removable.
Of course, the theorem implies that there are no such spaces.

\begin{lemma}
\label{lemma-nowhere}
Assume that $X$, $Y$, $n\ge 2$, $\pi$ and the
$K^j_y$ and $U^j_y$ are as in 
Definition \ref{def-superdup}, with  $X$ nowhere
$(n-1)$--\sdup dissipated and nowhere removable.
\begin{itemizz}
\item[1.] For a fixed $j$ and non-empty open $V \subseteq Y$:
$U^j_y \cap V \ne \emptyset$ for some $y \in V$.
\item[2.] For any $\varepsilon > 0$, the sets:
\begin{align*}
A^j_\varepsilon := \{z \in Y : \exists y [ z \in U^j_y
\ \&\ d(y,z) < \varepsilon ] \} \\
B^j_\varepsilon := \{y \in Y : \exists z [ z \in U^j_y
\ \&\ d(y,z) < \varepsilon ] \} 
\end{align*}
are dense and open in $Y$.
\end{itemizz}
\end{lemma}
\begin{proof}
For (1): Assume that  $U^j_y \cap V = \emptyset$ for all $y \in V$.
Let $W$ be a non-empty clopen subset of $\pi\iv\{V\}$, 
and consider the restriction 
$\pi\res W : W \onto \widetilde Y = \pi( W )$.
$\widetilde U^j_y = U^j_y \cap \widetilde Y = \emptyset$ for each
$y \in \widetilde Y$ and
$\widetilde K^j_y = K^j_y \cap W = ( K^j_y \setminus \pi\iv(U^j_y)  ) \cap W$
is empty for each $y \in \widetilde Y$
because it is clopen in $X$ and removable.
But then, by deleting index $j$, we see that
$W$ is $(n-1)$--\sdup dissipated;
in the special case $n = 2$, $W$ would be countable because $X$
does not contain a Cantor subset.

For (2):
They are open by (3) of Definition \ref{def-superdup}.
If one of them fails to be dense, then there is a non-empty open
$V \subseteq Y$ such that $V$ is disjoint from either
$A^j_\varepsilon$ or  $B^j_\varepsilon$.  In either case,
we may assume that $\diam(V) < \varepsilon$ which implies that
$z \notin U^j_y$ whenever $z,y \in V$, contradicting (1).
\end{proof}

\begin{lemma}
\label{lemma-super-big}
If $n \ge 2$ and $X$ is $n$--\sdup dissipated and 
nowhere $(n-1)$--\sdup dissipated and nowhere
removable, then $X$ is $n$--big.
\end{lemma}
\begin{proof}
Fix $\AA \lea C(X)$ and $\Upsilon : \AA \to \omega$.  We shall
verify the conclusion of Definition \ref{def-big} with $r = 1$,
so we shall find
$f_0, \ldots, f_{n-1} \in \AA$ and $c \in X$
such that
the $\Upsilon(f_j)$, for $j = 0, \ldots, n-1$, are all equal,
and such that $|f_i(c) - f_j(c)| \ge 1$ whenever $0 \le i < j < n$.
Let $Y$, $\pi$ and the
$K^j_y$ and $U^j_y$ be as in Definition \ref{def-superdup}.
Let $G = \bigcap\{A^j_\varepsilon \cap B^j_\varepsilon: 
\varepsilon > 0 \ \&\ j < n\}$; by Lemma \ref{lemma-nowhere},
$Y \backslash G$ is of first category in $Y$ because the intersection
may be taken just over rational $\varepsilon$.

If $y \in G$, then $y$ is in the closure of each $U^j_y$,
so that $\pi\iv\{y\}$ meets each $K^j_y$; let $x^j_y$ be
the element of $\pi\iv\{y\} \cap K^j_y$.
Since finite spaces have the CSWP, we may choose,
for each $y \in G$, a $g_y \in \AA$ such that 
$g_y(x^j_y) = 2j$ for each $j < n$.
Then, chose a rational $\varepsilon_y > 0$ such that
$|g_y(x) - 2j| < 1/2$ whenever $j < n$,
$x \in K^j_y$, and $d(\pi(x),y) < \varepsilon_y$.

Now, fix $N \subseteq G$, $\varepsilon > 0$, and $\ell \in \omega$
such that $N$ is not of first category in $Y$ and
$\varepsilon_y = \varepsilon$ and $\Upsilon(g_y) = \ell$ for all $y \in N$.
Then, fix a point $d \in N$ and a $\delta$ with
$0 < \delta < \varepsilon$ such that
$N \cap B(d, \delta)$ is dense in $B(d, \delta)$.
Let $c$ be any point in $\pi\iv\{d\}$.
For each $j < n$, $\{y : d \in U^j_y\}$ is open, and this
set meets $B(d, \delta)$ (since $d \in N \subseteq A^j_\delta$),
so choose $y^j \in N \cap B(d, \delta)$ such that
$d \in U^j_{y^j}$, and we can let $f_j = g_{y_j}$;
note that $d \in U^j_{y^j} \to c \in K^j_{y^j} \to |f_j(c) - 2j| < 1/2$.
\end{proof}

\begin{proofof}{Theorem \ref{thm-superdup-remov}}
Apply Theorem \ref{thm-R-remov};
every $n$--\sdup dissipated space $X$ is in $\RRRR$ by
Lemmas \ref{lemma-get-nowhere} and \ref{lemma-super-big}.
\end{proofof}

\section{Remarks and Questions}
\label{sec-rem}

Regarding our notion of bigness:
From the point of view of general topology,
the use of the  ``$\AA \lea C(X)$'' in Definition \ref{def-big}
seems a bit artificial, although it was needed for the CSWP proofs.
It would be more natural to restrict $\AA$ to be only $C(X)$,
which would result in a weaker property; but we do not
know if it would really be strictly weaker.
Of course, we can always replace $\AA$ by $\cl(\AA)$,
so the two properties are equivalent when $X$ has the CSWP.

The degree of bigness of some LOTSes is easily calculated.
Doing so lets us show (Example \ref{ex-big-limprop}) 
that the ``obvious'' generalization of
Lemma \ref{lemma-three-big} is false.
It is easy to see that $\omega_1 + 1$ is $n$--big for all $n$.
But there is a class of 
LOTSes for which the bigness is bounded.  We do not state the most general
possible result,  but just say enough to
verify Example \ref{ex-big-limprop}, which uses the $I_\Lambda$
from Definition \ref{def-das}.

\begin{lemma}
\label{lemma-small-big}
Let $L = I_\Lambda$, where $\Lambda : I \to \omega$, and let
$K$ be any compact space which is not $(n+1)$--big.  Let $X = L \times K$.
Then $X$ is not $(3n+1)$--big.
\end{lemma}
\begin{proof}
Let $\sigma: L \onto I$ be the standard map.
Also, applying the definition of ``not $(n+1)$--big'',
fix $\AA \lea C(K)$ and $\Upsilon : \AA \to \omega$
such that for each $c \in K$ and each
$f_0, f_1, \ldots, f_n \in \AA$  with
$\Upsilon(f_0) = \Upsilon(f_1) = \cdots = \Upsilon(f_n)$,
there are $j < k \le n$ such that $|f_j(c) - f_k(c)| < 1/4$.

Let $M = C(K)$, with the usual $\sup$ norm.
For $f \in C(X)$, define $\widetilde f \in C(L, M)$
by $(\widetilde f(u))(z) = f(u,z)$.
Let $\BB$ be the set of all $f \in C(X)$ such that
$\widetilde f (u) \in \AA$ for all $u \in L$.
Then $\BB \lea C(X)$, and we shall define a partition $\Phi$
of $\BB$ into $\aleph_0$ pieces demonstrating that $X$ is not $(3n+1)$--big.
As a first approximation,
for each $f \in \BB$, choose
$\Psi(f) = (m^f,  \vec y^f, \vec r^f, \vec s^f, \vec t^f)$ so that:
\begin{itemizz}
\item[1.] $1 \le m^f \in \omega$. 
\item[2.] $\vec y ^f = \langle y^f_i : 0 \le i \le 2m^f \rangle$,
each $y^f_i \in I$, and $y^f_i \in \QQQ$ when $i$ is even.
\item[3.] $0 = y^f_0 < y^f_1 < \cdots < y^f_{2m^f} = 1$.
\item[4.] $\vec r ^f = \langle r^f_i : 0 \le i \le 2m^f \rangle$,
where each $r^f_i = |\sigma\iv\{y^f_i\}| - 1 = \Lambda(y^f_i)$.
\item[5.] $\| \widetilde f(u) - \widetilde f(v)\| \le 1/4$ whenever
$ \max(\sigma\iv\{y^f_i\}) \le u \le v \le \min(\sigma\iv\{y^f_{i+1}\})$.
\item[6.] $\vec s ^f = \langle s^f_{i,\mu} :
0 \le i \le 2m^f \ \&\  0 \le \mu \le r^f_i\rangle$, where
$\{s^f_{i,\mu} : 0 \le \mu \le r^f_i\} \subset L$
lists $\sigma\iv\{y^f_i\}$ in increasing order;
so $s^f_{i,\mu} = (y^f_i, \mu)$.
\item[7.] $\vec t ^f = \langle t^f_{i,\mu} :
0 \le i \le 2m^f \ \&\  0 \le \mu \le r^f_i\rangle$, where
$t^f_{i,\mu} = \Upsilon(\widetilde f ( s^f_{i,\mu} ))$.
\end{itemizz}
Such a $\Psi(f)$ may be chosen using compactness, plus  continuity
of $\widetilde f$.
Of course, there are $2^{\aleph_0}$ possible values of $\Psi(f)$
because of the $y^f_i$ and $s^f_{i,\mu}$ for odd $i$, so we
delete these and define
$\Phi(f) = (m^f,  \vec \yyy^f, \vec r^f, \vec \sss^f, \vec t^f)$, where

\begin{itemizz}
\item[8.] $\vec \yyy ^f = \langle y^f_i : 0 \le i \le 2m^f 
\ \&\ \mbox{$i$ is even} \rangle$.
\item[9.] $\vec \sss ^f = \langle s^f_{i,\mu} :
0 \le i \le 2m^f \ \&\   \mbox{$i$ is even}\ \& \ 1 \le \mu \le r^f_i \rangle$.
\end{itemizz}

\noindent
There are only countably many possible values for $\Phi(f)$,
so if $X$ were $(3n+1)$--big, we could fix a $(b,c) \in X = L \times K$
and $f_0, \ldots, f_{3n} \in \AA$  such that
the $\Phi(f_j)$  are all the same, and such that 
$|f_j(b,c) - f_k(b,c)| \ge 1$ whenever $j < k \le n$.
We shall now derive a contradiction.
Write $\Phi(f_j) = (m,  \vec \yyy, \vec r, \vec \sss, \vec t)$.

If $b =  s^f_{i,\mu}$ for some even $i$, then
the $\Upsilon(\widetilde f_j (b)) = t^f_{i,\mu}$ are all the same, and 
we contradict our assumptions on $\Upsilon$ just using
$\widetilde f_0(b) , \ldots , \widetilde f_n(b)$.
So, we may fix an even $i < 2m$ so that
$ \max(\sigma\iv\{y_i\}) <  b <  \min(\sigma\iv\{y_{i+2}\})$.
Now, for each $j \in \{0, 1, \ldots, 3n\}$, there are three cases:
\begin{itemizz}
\item[I.]$ \max(\sigma\iv\{y_i\}) <  b <  \min(\sigma\iv\{y^{f_j}_{i+1}\})$.
\item[II.] $b \in \sigma\iv\{y^{f_j}_{i+1}\}$.
\item[III.]$ \max(\sigma\iv\{y^{f_j}_{i+1}\}) <  b <  
    \min(\sigma\iv\{y_{i+2}\})$.
\end{itemizz}
So, one of these cases must happen for $n+1$ values of $j$.
We shall assume that this is Case I, since the argument is essentially the
same in the other two cases.
Permuting the $f_j$, we may assume that Case I holds for $0 \le j \le n$.
Fix $\mu = r_i$, so that $\max(\sigma\iv\{y_i\}) = s_{i,\mu}$,
so $\Upsilon(\widetilde f_j (s_{i,\mu})) = t^f_{i,\mu}$ for each $j \le n$.
By our assumptions on $\Upsilon$, we may fix
$j < k \le n$ such that
$| f_j (s_{i,\mu}, c) - f_k (s_{i,\mu}, c) | < 1/4$.
Applying Condition (5) above, we have
$| f_j (b, c) - f_k (b, c) | < 3/4$,
contradicting $|f_j(b,c) - f_k(b,c)| \ge 1$.
\end{proof}

In particular, letting $K$ be the $1$--point space, we see
that an  $I_\Lambda$ is not $4$--big.  Then, proceeding by induction,

\begin{lemma}
\label{lemma-small-lotsprod}
$\prod_{j<m} L_{\Lambda_j}$ is not $(3^m + 1)$--big.
\end{lemma}

We remark that in the proof of Lemma \ref{lemma-small-big},
we could have replaced the ``$<$'' by ``$\le$'' in Cases I and III,
although then they would not be disjoint from Case II.
However, in the special case of $L = I_S$, Case II can now
be eliminated, so that we can replace the
``$(3n+1)$--big'' by ``$(2n+1)$--big'', obtaining:

\begin{lemma}
\label{lemma-smaller-big}
Let $L = I_S$, where $S \subseteq I$, and let
$K$ be any compact space which is not $(n+1)$--big.  Let $X = L \times K$.
Then $X$ is not $(2n+1)$--big.
\end{lemma}

\begin{lemma}
\label{lemma-das-not-big}
$\prod_{j<m} I_{S_j}$ is not $(2^m + 1)$--big.
\end{lemma}

\begin{example}
\label{ex-big-limprop}
For any $n > 3$, there is an $X$ which 
does not have a Cantor subset and which is not $7$--big,
such that $X$ is $(n + 2)$--superdissipated and not 
$(n + 1)$--superdissipated.
\end{example}
\begin{proof}
For $n \ge 1$, let $L_n = I_{\Lambda_n}$, where
$\Lambda_n(x) = n$ for $n \in (0,1)$, and 
$\Lambda_n(0) = \Lambda_n(1) = 0$.
Then $L_1$ is the double arrow space.
Let $X = L_n \times L_1$.  Then $X$ is not $7$--big
by Lemmas \ref{lemma-small-lotsprod} and \ref{lemma-smaller-big}.
$X$ is  $(n + 2)$--dissipated by Lemma \ref{lemma-tight-prod}.
To prove that $X$ is not $(n + 1)$--dissipated, it is sufficient
(by Lemma 3.6 of \cite{KU2}) to observe that 
for each $\varphi \in C(L_n, [0,1]^\omega)$
there is a $z \in [0,1]^\omega$ with $|\varphi\iv\{z\}| \ge n+1$.
\end{proof}

It is easily seen using Lemma \ref{lemma-big-helper}
that $(I_S)^n$ is $2^n$--big when $S$
is uncountable, so it has the NTIP, since it is also
$(2^{n-1} + 1)$--superdissipated.  However, it is not clear whether
it has the CSWP in the case that $S$ meets all Cantor sets, since the
natural proof requires looking at arbitrary 
perfect subspaces of $(I_S)^n$.

\begin{example}
\label{ex-seven}
If $S \subseteq (0,1)$ meets all Cantor sets, then
there is a perfect $X \subseteq (I_S)^4$
such that $X$ is not $8$--dissipated, is not $7$--big,
and has no Cantor subsets.
\end{example}
\begin{proof}
Let $D \subset (I_S)^3$ be the diagonal, and let $X = D \times I_S$.
Then $D$ is the same as the LOTS obtained from $I$ by replacing each point
in $S$ by eight points.  Since $7 = 2\cdot 3 + 1$,
Lemmas \ref{lemma-small-lotsprod} and \ref{lemma-smaller-big} show
that $X$ is not $7$--big.
The proof of Example \ref{ex-big-limprop} shows that
$X$ is not $8$--dissipated.
\end{proof}

This particular $X$ has the CSWP, and in fact is removable
(see Section \ref{sec-remov}), since $I_S$ is removable,
and after removing the clopen copies of $I_S$ from $X$,
we are left with a copy of $(I_S)^2$, which is also removable.
We do not know whether $(I_S)^4$ itself must have the CSWP.

A simple example of $n$--big spaces is given by:

\begin{proposition}
If $X$ is compact and $|X| > 2^{\aleph_0}$, then $X$ is $n$--big for 
all $n\in\omega$.
\end{proposition}
\begin{proof}
Fix $\AA \lea C(X)$, fix $n \in \omega$, and fix
$\Upsilon : \AA \to \omega$.  We shall verify the conclusion
of Definition \ref{def-big} with $r = 1$.

Let $\PPP$ be the set of all finite partial functions from 
$X$ to $\omega$; so each $p \in \PPP$ is a function with  $\dom(p)$ a
finite subset of $X$ and $\ran(p) \subseteq \omega$.
For $p \in \PPP$, choose an $f_p \in \AA$ with $p \subset f_p$.

For each $c \in X$ and $s \in \omega$, 
let $E_{c,s} = \{p(c) :
p \in \PPP \ \&\ c \in \dom(p) \ \&\ \Upsilon(f_p)
= s\}\allowbreak\subseteq \omega$.
If some $|E_{c,s}| \ge n$ then we are done, so assume
that $|E_{c,s}| \le n-1$ for all $c,s$.
There are only $ 2^{\aleph_0}$ possibilities for
$\langle E_{c,s} : s \in \omega \rangle$, so we can fix an infinite
$A \subseteq X$
and sets $E_s \in [\omega]^{<n}$ for $s\in \omega$ such that 
$E_{c,s} = E_s$ for all $c \in A$ and all $s \in \omega$.
But then $\ran(p) \subseteq E_s$ whenever $p \in \PPP$ and
$\Upsilon(f_p) = s$ and $\dom(p) \subseteq A$.
Now choosing $p$ with $\dom(p) \subseteq A$ and 
$|\ran(p)| = n$ yields a contradiction.
\end{proof}

Finally, the following Ramsey-type lemma might be of interest
for studying products of LOTSes, although we never needed it
in this paper.  The proof uses the terminology from
Definition \ref{def-boxes}.

\begin{proposition}
Fix an uncountable $J \subseteq \RRR^n$, and assume that
\[
\forall x,y \in J \, [ x \ne y \to \forall i < n [x_i \ne y_i]]
\tag{$*$} \ \ .
\]
Then there is a 1-1 function $\varphi: \QQQ \to J$ such that for all
$i < n$:
\[
\forall p,q \, [ p < q \to \varphi(p)_i < \varphi(q)_i ]
\mbox{\ \ {\rm or} \ \ }
\forall p,q \, [ p < q \to \varphi(p)_i > \varphi(q)_i ] \ \ .
\tag{\dag}
\]
\end{proposition}
\begin{proof}
Call a box $B  = \bx[a,b]$ 
\emph{big} iff $B \cap J$ is uncountable;
By  $(*)$, this implies that
$B^\circ \cap J$ is uncountable, where $B^\circ$ denotes
the interior of $B$.
For $\Delta \in \DD$, let $-\Delta$ result from interchanging
the signs $+$ and $-$ in $\Delta$.
Call the box $B$ $\Delta$--\emph{bad} iff $B$ is big and
there is no $d \in B^\circ \cap J$ 
such that $\corn(d, B, \Delta)$ and $\corn(d, B, -\Delta)$ 
are both big.  Observe, for any big box $B$:

\begin{itemizz}
\item[1.] $B$ is $\Delta$--bad iff $B$ is $(-\Delta)$--bad. 
\item[2.] If $B$ is $\Delta$--bad and $A \subseteq B$ is a big box,
then $A$ is $\Delta$--bad. 
\item[3.] There is some $\Delta\in\DD$ such that $B$ is not $\Delta$--bad.
\end{itemizz}
(1) and (2) are obvious.  To prove (3),  we note first that if we
replaced $\QQQ$ by $\omega$ or a finite set in the statement of the lemma,
then the
result would be obvious by Ramsey's Theorem. 
Now, let $Z$ be the set of
points of $B^\circ \cap J$ which are condensation points of $J$.
Obtain $\varphi: \{0,1,2\} \to Z$ so that (\dag) holds replacing
$\QQQ$ by $\{0,1,2\}$.
Let $a = \varphi(0)$, $b = \varphi(2)$, and $d = \varphi(1)$.
By (\dag), there is some $\Delta \in \DD$ such that
$a \in \corn(d, B, \Delta)$ and $b \in \corn(d, B, -\Delta)$,
and then $\corn(d, B, \Delta)$ and $\corn(d, B, -\Delta)$ are both big.

Using (2) (sub-boxes go from bad to worse) and (3), we can
fix a $\Delta \in \DD$ and a big box $B$ such that
for all big boxes $A \subseteq B$,
$A$ is not $\Delta$--bad.    We may now list $\QQQ$ in type
$\omega$ and obtain $\varphi$ in $\omega$ steps.
When $p < q$,
we shall have $\varphi(p)_i < \varphi(q)_i $ when $\Delta_i = +1$ and
$\varphi(p)_i > \varphi(q)_i $ when $\Delta_i = -1$.
\end{proof}

\end{document}